# STANDARD DEVIATION OF THE LONGEST COMMON SUBSEQUENCE

By Jüri Lember[1] and Heinrich Matzinger

*University of Tartu and Georgia Tech*

Let $L_n$ be the length of the longest common subsequence of two independent i.i.d. sequences of Bernoulli variables of length $n$. We prove that the order of the standard deviation of $L_n$ is $\sqrt{n}$, provided the parameter of the Bernoulli variables is small enough. This validates Waterman's conjecture in this situation [*Philos. Trans. R. Soc. Lond. Ser. B* **344** (1994) 383–390]. The order conjectured by Chvatal and Sankoff [*J. Appl. Probab.* **12** (1975) 306–315], however, is different.

**1. Introduction.** Throughout this paper $X_1, X_2, \ldots$ and $Y_1, Y_2, \ldots$ are two independent sequence of i.i.d. Bernoulli variables with parameter $0.5 \geq \varepsilon > 0$:

$$\varepsilon = P(X_i = 1) = P(Y_i = 1) = 1 - P(X_i = 0) = 1 - P(Y_i = 0).$$

Let $X := X_1 X_2 \cdots X_n$ and let $Y := Y_1 Y_2 \cdots Y_n$. The longest common subsequence (LCS) of $X$ and $Y$ is any common subsequence that has the longest possible length. The length of LCS is denoted $L_n$. Formally, $L_n$ is the biggest $k$ such that there exists two subsets of indices $\{i_1, \ldots, i_k\}, \{j_1, \ldots, j_k\} \subset \{1, \ldots, n\}$ satisfying $i_1 < i_2 < \cdots < i_k$, $j_1 < j_2 < \cdots < j_k$ and $X_{i_1} = Y_{i_1}, X_{i_2} = Y_{i_2}, \ldots, X_{i_k} = Y_{i_k}$. The main result of this paper is, that for $\varepsilon > 0$ small enough, the order of the standard deviation of $L_n$ is $\sqrt{n}$.

LCS's are a very important tool in computational biology, where they are used for comparing DNA- and protein-alignments (see, e.g., [3, 16, 17]). They are also used in computational linguistics, speech recognition and so on. In all these applications, two strings with a relatively long LCS, are deemed related.

---

Received April 2006; revised February 2007.

[1]Supported by Estonian Science Foundation Grant 7553, and SFB701 of Bielefeld University

EXAMPLE. Let us give an example of the practical use of LCS's. Take the two words: $X = fanthastic$ and $Y = fntastique$. These two words are very similar. They were obtained from the English word "fantastic" and the French word "fantastique" by adding spelling mistakes. We would like the computer to recognize the similarity. If the computer compares letter by letter,

| f | a | n | t | h | a | s | t | h | a | s | t | i | c |
|---|---|---|---|---|---|---|---|---|---|---|---|---|---|
| f | n | t | a | s | t | i | q | u | e |   |   |   |   |

,

it finds that only one letter coincides. Comparing the $i$th letter of the first word with the $i$th letter of the second word for all the letters is not a good way to recognize any similarity. The reason are the missing letters. The original position of the letters in the words gets changed. To take into account the missing letters or added letters, we align the two words allowing for gaps. We allow only same letters to be matched with each other. In such a way, we obtain a sequence of letters that is contained in $X$ as well as in $Y$. Such a sequence is a common subsequence of $X$ and $Y$. Hence, the longest common subsequence is the maximum number of same letters we can align allowing gaps. In our example the maximum is given by the alignment

(1)

| f | a | n | t | h | a | s | t | i | c |   |   |   |
|---|---|---|---|---|---|---|---|---|---|---|---|---|
| f |   | n | t |   | a | s | t | i |   | q | u | e |

.

Hence $f, n, t, a, s, t, i$ is the longest common subsequence of the two words and the length of the longest common subsequence, $L_n$, is 7. This indicates that the two words are very similar.

To distinguish related pairs of strings from unrelated via the LCS method, we need to assess the order of the fluctuation of the LCS. For this reason the random variable $L_n$ has received a lot of attention. Nonetheless, many questions remain open. In their pioneering paper [7], Chvatal and Sankoff prove that the limit

$$\gamma := \lim_{n \to \infty} \frac{EL_n}{n} \qquad (2)$$

exists. In [1], Alexander investigated the rate of the convergence in (2) and showed that for a constant $C$, $EL_n - n\gamma \geq C\sqrt{n \ln n}$. Moreover, by a subadditivity argument

$$\frac{L_n}{n} \to \gamma \qquad \text{a.s and in } L_1 \qquad (3)$$

(see, e.g., [1, 17]). The constant $\gamma$ is called the Chvatal–Sankoff constant and its value is unknown for even as simple cases as i.i.d. Bernoulli sequences. In this case, the value of $\gamma$ obviously depends on the Bernoulli parameter

$\varepsilon$. When $\varepsilon = 0.5$, the various bounds indicate that $\gamma \approx 0.81$ [4, 11, 14]. For a smaller $\varepsilon$, $\gamma$ is even bigger. Further bounds on $\gamma$ have been obtained by Martinez, Hauser and Matzinger [9]. Hence, a common subsequence of two independent Bernoulli sequences typically makes up a large part of the total length. This implies that to make some inference, the size of the variance $\text{Var}[L_n]$ is essential. Unfortunately, not much is known about $\text{Var}[L_n]$ and its asymptotic order is one of the central open problems in string matching theory. Monte Carlo simulations lead Chvatal and Sankoff in [7] to conjecture for $\varepsilon = 0.5$ that $\text{Var}[L_n] = o(n^{2/3})$. Using an Efron–Stein type of inequality, Steele [14] proved $\text{Var}[L_n] \leq 2\varepsilon(1-\varepsilon)n$. In [15], Waterman asks whether this linear bound can be improved. His simulations suggest that for $\varepsilon < \frac{1}{2}$ this is not the case and $\text{Var}(L_n)$ grows linearly. In [6], Boutet de Monvel simulates $\text{Var}(L_n)$ for the case $\varepsilon = \frac{1}{2}$ and notices the linear growth as well. However, he adds that the linear regime of the growth is not reached before $n$ is about 10,000. He also simulates the values of the random variable

$$\frac{L_n - EL_n}{\sqrt{\text{Var}(L_n)}}$$

and founds its distribution close to normal.

In a series of papers, we investigate the asymptotic behavior of $\text{Var}[L_n]$ in various setups. Our goal is to find out, whether there exists a constant $c > 0$ (not depending on $n$) such that $\text{Var}[L_n] \geq cn$. Together with Steele's bound, this means that $cn \leq \text{Var}[L_n] \leq n$, that is, $\text{Var}[L_n] = \Theta(n)$ [a sequence $a_n$ is of order $\Theta(n)$, if, for some constants $0 < c < C < \infty$, $cn \leq a_n \leq Cn$ for all $n$ large enough]. In [5], Bonetto and Matzinger consider the asymmetric situation where the random variables in $X$ are Bernoulli with $1/2$, but $Y$ is a random i.i.d. string with three symbols. They obtain that in this setting $\text{Var}[L_n] = \Theta(n)$. In [10], Houdre, Lember and Matzinger investigate the asymptotic behavior of the longest common increasing subsequence of two independent Bernoulli sequences (a binary increasing sequence begins with a block of zero's followed by a block of one's). They find that under this additional restriction $n^{-1/2}(L_n - EL_n)$ converges in law to a functional of two Brownian motions implying that $\text{Var}[L_n] = \Theta(n)$ holds again (here $L_n$ designates the length of the longest common increasing subsequence). Durringer, Lember and Matzinger [13] show that $\text{Var}[L_n] = \Theta(n)$ when $Y$ is a nonrandom periodic binary sequence and $X$ is an i.i.d. Bernoulli $1/2$ sequence. The nature of the optimal path has been investigated by Amsalu, Popov and Matzinger in [2] as well as by Lember, Matzinger and Vollmer in [12].

The relatively long history shows that determining the exact order of the fluctuation of $L_n$ is a difficult problem. In fact, as noted in [1, 3], the LCS problem can be reformulated as a Last Passage Percolation (LPP) problem with correlated weights. But for standard LPP and First Passage

Percolation, the question of the exact order of the fluctuation remain open except for the case of geometric or exponential weights which has been solved by Johanson.

**2. Main result.** The main result of this paper, Theorem 2.1, asserts that when $\varepsilon > 0$ is small, the fluctuation of $L_n$ is of order $\sqrt{n}$. In fact, the theorem gives only a lower linear bound for the variance of $L_n$. The upper linear bound comes from the result of Steele [14]. Hence, Theorem 2.1 implies that $\mathrm{Var}[L_n] = \Theta(n)$.

THEOREM 2.1.   *There exists $\varepsilon_0 > 0$ such that for every $\varepsilon < \varepsilon_0$, there exists a constant $c > 0$ depending on $\varepsilon$ but not depending on $n$, that satisfies*

$$\mathrm{Var}[L_n] \geq c \cdot n \qquad \forall n.$$

One of the main tools in this paper is a map that picks a one in the text $X$ or $Y$ at random and changes it into a zero. Let $\tilde{X}$ and $\tilde{Y}$ designate the texts obtained in this way.

EXAMPLE.   Let $n = 6$, $X = 001000$ and $Y = 101000$. The total number of ones in the two texts is 3. Hence, we pick one of these three ones at random with equal probability and switch it into a zero. Assume we pick the second one in text $Y$. Then $\tilde{X} = 001000$ and $\tilde{Y} = 100000$.

Let us define $\tilde{X}$ and $\tilde{Y}$ rigorously. For a binary string $x = x_1 x_2 \cdots x_n$, we denote by $N_1^x$ the total number of ones in $x$. So $N_1^x := \sum_{i=1}^n x_i$. Similarly, $N_1^y$ is the total number of ones in $y = y_1 y_2 \cdots y_n$. The binary random strings $\tilde{X}$ and $\tilde{Y}$ are defined by the following equations:

$$\sum_{i=1}^n (|\tilde{X}_i - X_i| + |\tilde{Y}_i - Y_i|) = \begin{cases} 1, & \text{if } \sum_{i=1}^n (X_i + Y_i) > 0; \\ 0, & \text{else,} \end{cases}$$

$$\sum_{i=1}^n (\tilde{X}_i - X_i + \tilde{Y}_i - Y_i) = \begin{cases} -1, & \text{if } \sum_{i=1}^n (X_i + Y_i) > 0; \\ 0, & \text{else,} \end{cases}$$

$$P(\tilde{X}_i \neq X_i | X = x, Y = y) = \begin{cases} 0 & \text{if } x_i = 0; \\ \dfrac{1}{N_1^x + N_1^y}, & \text{else,} \end{cases}$$

$$P(\tilde{Y}_i \neq Y_i | X = x, Y = y) = \begin{cases} 0 & \text{if } y_i = 0; \\ \dfrac{1}{N_1^x + N_1^y}, & \text{else.} \end{cases}$$

Let $\tilde{L}_n$ denote length of the longest common subsequence of $\tilde{X}$ and $\tilde{Y}$. When we change one bit in $X$ or $Y$ and flip it to the opposite value, then the length of the LCS changes by at most one. The next theorem shows that in this case the length of the LCS $L_n$ is more likely to increase by one unit than to decrease by one unit.

THEOREM 2.2. *There exist constants $\alpha_1$ and $\alpha_2$, $\alpha_1 > \alpha_2$ and a set $B_n \subset \{0,1\}^n \times \{0,1\}^n$ such that for all $(x,y) \in B_n$*

(4) $$P(\tilde{L} - L = 1 | X = x, Y = y) \geq \alpha_1,$$

(5) $$P(\tilde{L} - L = -1 | X = x, Y = y) \leq \alpha_2.$$

*Moreover, there exists an $\varepsilon_0 > 0$ such that for every $0 < \varepsilon \leq \varepsilon_o$*

(6) $$P((X,Y) \in B_n) \geq 1 - e^{-c_1 n},$$

*where $c_1 > 0$ does not depend on $n$, but may depend on $\varepsilon$.*

In Section 3, we prove that Theorem 2.2 implies Theorem 2.1. Let us briefly explain the main ideas behind the proof. We define two sequences of random binary strings $X^1, X^2, \ldots, X^{2n}$ and $Y^1, Y^2, \ldots, Y^{2n}$, all of them having length $n$. The strings $X^k$ and $Y^k$ are define by induction on $k$: $X^{2n}$ and $Y^{2n}$ consist only of ones; $X^{k-1}$ and $Y^{k-1}$ are obtained by choosing a one at random in $X^k Y^k$ and replacing it by a zero. Hence we use the random map $\tilde{\ }$. We designate by $L(k)$ the length of the LCS of $X^k$ and $Y^k$. Note that the total number of ones in the string $X^k$ and $Y^k$ is $k$. Let $(X,Y)$ be independent of $\{(X^k, Y^k)\}_{k \in \{0,\ldots,2n\}}$ and let $N_1$ designate the total number of ones in the two strings $X$ and $Y$. It is not hard to see that $(X^k, Y^k)$ has the same distribution as $(X,Y)$ conditional on $N_1 = k$. This implies that $L(N_1)$ has same distribution, as $L_n$. The standard deviation of $N_1$ is of order $\sqrt{n}$. Moreover, from Theorem 2.2 directly follows that the (random) map $k \mapsto L(k)$ tends to increase linearly on a certain scale. These two facts together imply immediately that the standard deviation of $L(N_1)$ and hence also of $L_n$ is of order $\sqrt{n}$.

Let us now give a heuristic argument why Theorem 2.2 holds. Recall that in this paper, we consider the situation where one has a small, but fixed probability. Hence, in the texts $X$ and $Y$, there is a small proportions of ones. This implies that only a small percentage of ones can figure in a LCS. It will turn out that the number of ones in a LCS is typically of order $\varepsilon^2 n$. This is much less than the total number of ones in the texts $X$ and $Y$, which is of order $2\varepsilon n$. It follows that the majority of ones in the texts $X$ and $Y$ constitute a "net loss" for the score $L_n$. Hence the number of ones tends to influence the score $L_n$ negatively. Changing a randomly picked one into zero is not very likely to decrease the score. It can decrease the score only if the chosen one is used in a LCS. But the additional zero obtained in this way will in many cases increase the score.

EXAMPLE. Let

$$X = 0001000010000000000001,$$

$$Y = 00010000000010000100000.$$

The longest common subsequence $Z$ is $Z = 000100000000000000000$. An alignment corresponding to $Z$ is

| X | 0 | 0 | 0 | 1 | 0 | 0 | 0 | 0 | 1 | 0 | 0 | 0 | 0 |   | 0 | 0 | 0 | 0 |   | 0 | 0 | 0 | 0 | 0 | 1 |
|---|---|---|---|---|---|---|---|---|---|---|---|---|---|---|---|---|---|---|---|---|---|---|---|---|---|
| Y | 0 | 0 | 0 | 1 | 0 | 0 | 0 | 0 |   | 0 | 0 | 0 | 0 | 1 | 0 | 0 | 0 | 0 | 1 | 0 | 0 | 0 | 0 | 0 |   |
| Z | 0 | 0 | 0 | 1 | 0 | 0 | 0 | 0 |   | 0 | 0 | 0 | 0 |   | 0 | 0 | 0 | 0 |   | 0 | 0 | 0 | 0 | 0 |   |

The optimal solution is obtained by matching all the zeros, and the first one in both texts, but discarding all other ones. We see the general phenomena: since there are few ones, sometimes by chance some ones appear in respective positions in the two texts where they can be matched. The other ones in text $X$ and $Y$ appear in places in the text where we cannot match them with a one. If we would match them we would loose too many zeros. That is why, most ones can not be used in the LCS.

The argument in the previous numerical example gives a first idea of what is happening. However, proving anything rigorously is difficult. The reason is as follows. We take $\varepsilon$ small but fixed and let then $n$ tend to infinity. The optimal alignment (optimal alignment is the alignment which defines the LCS) is then going to be a global alignment. This means that typically some parts of the text $X$ will be connected with parts of the text $Y$ that are far away. This introduces complicated correlations between the different parts of the optimal alignment. Microscopically it is easy to understand the approximate behavior of the optimal alignment. Macroscopically however, little is understood about the optimal alignment. It seems that there are complicated long range interactions between all the different parts.

**3. Theorem 2.2 implies Theorem 2.1. The proof.** In this section, we prove that Theorem 2.2 implies Theorem 2.1. We use some of the techniques developed in [5].

Recall that $N_1$ is the total number of ones in the two strings $X$ and $Y$. We already mentioned briefly the definition of the random pair of strings $(X^k, Y^k)$ for $k \in [0, 2n]$. Let us give more details. Both strings $X^k$ and $Y^k$ are binary strings of length $n$. We proceed recursively on $k$. The strings $X^{2n}$ and $Y^{2n}$ consist only of 1's. We pick a 1 in the strings $X^{2n}Y^{2n}$ at random and change it into a 0. This way we obtain $(X^{2n-1}, Y^{2n-1})$. For general $k$, we obtain $(X^{k-1}, Y^{k-1})$ from $(X^k, Y^k)$ by choosing a 1 at random in $X^k Y^k$ and changing it to the opposite value. Each one has the same probability to get chosen. We request that conditional on $(X^k, Y^k)$, which one in $(X^k, Y^k)$ gets chosen, is independent of $\{(X^i, Y^i)\}_{i \in [k, 2n]}$. In other words, we apply the transformation $\tilde{\ }$, so that

$$X^{k-1} := \tilde{X}^k \quad \text{and} \quad Y^{k-1} := \tilde{Y}^k.$$

The distribution of $(X^k, Y^k)$ is equal to the distribution of $(X, Y)$ conditional on $N_1 = k$:

(7) $$\mathcal{L}(X^k, Y^k) = \mathcal{L}(X, Y | N_1 = k),$$

where $\mathcal{L}(W)$ designates the distribution of the random variable $W$.

Let $L(k)$ designate the length of the LCS of $X^k$ and $Y^k$.

We assume that $\{X^k, Y^k\}_{k \in [0, 2n]}$ are independent of the random variable $N_1$. Picking $N_1$ according to its distribution gives us random strings $(X^{N_1}, Y^{N_1})$ that have the same distribution as $(X, Y)$. Therefore, the length $L(N_1)$ of the LCS of $(X^{N_1}, Y^{N_1})$, has the same distribution as $L_n$. Hence

$$\text{Var}[L_n] = \text{Var}[L(N_1)].$$

Recall that our aim is to prove that $\text{Var}[L(N_1)]$ it at least of order $n$. This follows from two facts: (1) the order of $\text{Var}[N_1]$ is $n$; (2) the (random) map $k \mapsto L(k)$ typically decreases linearly on a certain scale.

The second point follows rather directly from Theorem 2.2 and is proven in Lemma 3.2. This section is dedicated to showing that (1) and (2) above imply the linear lower bound for $\text{Var}[L(N_1)]$. There are two technical difficulties: (a) the map $k \mapsto L(k)$ does not increase at every point, but only on a certain scale; (b) the increasing slope on a certain scale only holds in a domain where typically $N_1$ takes values, but not everywhere.

Recall that for any variables $V$ and $W$,

(8) $$\text{Var}[V] = \text{Var}[E[V|W]] + E[\text{Var}[V|W]] \geq E[\text{Var}[V|W]],$$

where $\text{Var}[V|W]$ is the variance of the conditional distribution $\mathcal{L}(V|W)$. Applying (8) to our case, we find

(9) $$\text{Var}[L(N_1)] \geq E[\text{Var}[L(N_1)|L(\cdot)]],$$

where $L(\cdot)$ is the (random) map $k \mapsto L(k)$. Note that $N_1$ is independent of $L(\cdot)$.

Let $I$ be the interval

(10) $$I := [2\varepsilon n - \sqrt{\varepsilon(1-\varepsilon)2n}, 2\varepsilon n + \sqrt{\varepsilon(1-\varepsilon)2n}].$$

Let $\tilde{N}_1$ be a random variable, independent of $L(\cdot)$ and having the distribution of $N_1$ conditioned on $N_1 \in I$. From (8), it follows for every fixed $L$ that

$$\text{Var}[L(N_1)] \geq \text{Var}[L(N_1)|N_1 \in I]P(N_1 \in I) = \text{Var}[L(\tilde{N}_1)]P(N_1 \in I).$$

Hence, since $L$ and $\tilde{N}_1$ are independent,

(11) $$E[\operatorname{Var}[L(N_1)|L(\cdot)]] \geq E[\operatorname{Var}[L(\tilde{N}_1)|L(\cdot)]]P(N_1 \in I).$$

Assume that $f : \mathbb{R} \to \mathbb{R}$ is map such that, for a constant $c > 0$, $f'(x) > c$ for all $x \in \mathbb{R}$. Then, for any random variable $Y$, we have

(12) $$\operatorname{Var}[f(Y)] \geq c^2 \operatorname{Var}[Y].$$

(See Lemma 3.2 in [5] for the proof.) Hence, if the map $L(\cdot)$ had positive slope everywhere larger than $c > 0$, it would follow that $\operatorname{Var}[L(N_1)] \geq c \cdot \operatorname{Var}[N_1]$. Typically, the (random) map $k \mapsto L(k)$ does not strictly increase for every $k \in [0, 2n]$. But it is likely that in $I$ it increases by a linear quantity. We are next going to formulate a lemma, proven in [5] (Lemma 3.3 in [5]), which is a modification of inequality (12), for when the map $k \mapsto f(k)$ does not increase every $k$, but has a tendency to increase on some scale.

LEMMA 3.1. *Let $c, m > 0$ be two constants. Let $f : I \to \mathbb{Z}$ be a non decreasing map that satisfies the following conditions*

(13)     $f(j) - f(i) \leq (j - i) \qquad \forall i < j,$

(14)     $f(j) - f(i) \geq c \cdot (j - i) \qquad \forall i, j \text{ such that } i + m \leq j.$

*Let $B$ be an $I$-valued random variable such that $E|f(B)| < \infty$. Then*

(15) $$\operatorname{Var}[f(B)] \geq c^2 \left(1 - \frac{2m}{c\sqrt{\operatorname{Var}[B]}}\right) \operatorname{Var}[B].$$

Recall the definition of $I$ in (10). Let $\alpha_1$ and $\alpha_2$ be the constants from Theorem 2.2 and let $E^n_{\text{slope}}$ denote the event that for all $i, j \in I$, such that $i + n^{0.1} \leq j$, we have

(16) $$L(j) - L(i) \geq \alpha_3|i - j|,$$

where

$$\alpha_3 := \frac{\alpha_1 - \alpha_2}{2}.$$

In other words, the event $E^n_{\text{slope}}$ says that $L(\cdot)$ has a slope of at least $\alpha_3$ on $I$, when we look only at points which are at least $n^{0.1}$ away from each other. The next lemma shows that the event $E^n_{\text{slope}}$ has high probability, provided Theorem 2.2 holds.

LEMMA 3.2. *For a constant $c_4 > 0$,*

(17) $$P(E^n_{\text{slope}}) \geq 1 - e^{c_4 \cdot n^{0.1}},$$

*provided $n$ is sufficiently big.*

PROOF. Let $A_n^k$ denote the event that the random vector $(X^k, Y^k)$ takes the values in the set $B_n$ defined in Theorem 2.2. So

$$A_n^k := \{(X^k, Y^k) \subset B_n\}.$$

Let $A_n^{\text{all}}$ be the event

$$A_n^{\text{all}} := \bigcap_{k \in I} A_n^k.$$

Let

$$\Delta^k := \begin{cases} L(k-1) - L(k), & \text{when } A_n^k \text{ holds;} \\ 1, & \text{else.} \end{cases}$$

Let $i < j$ and consider the random variable

$$\sum_{k=i+1}^{j} \Delta^k.$$

When $(X^k, Y^k) = (x, y) \in B_n^k$, that is, $A_n^k$ holds, then Theorem 2.2 says that

$$P(\Delta^k = 1 | X^k = x, Y^k = y) \geq \alpha_1,$$
$$P(\Delta^k = -1 | X^k = x, Y^k = y) \leq \alpha_2,$$

implying that $E[\Delta^k | A_n^k, X^k, Y^k] \geq \alpha_1 - \alpha_2$. Since $E[\Delta^k | (A_n^k)^c] = 1 > \alpha_1 - \alpha_2$, we get

(18) $$E(\Delta_k | X^k, Y^k) \geq \alpha_1 - \alpha_2.$$

Let, for every $k = 2n+1, \ldots, 2$,

$$\mathcal{F}_k := \sigma(X^{2n}, Y^{2n}, \ldots, X^{k-1}, Y^{k-1}).$$

These $\sigma$-algebras perform a (reversed) filtration, because

$$\mathcal{F}_{2n+1} \subset \mathcal{F}_{2n} \subset \cdots \subset \mathcal{F}_2.$$

The random variable $\Delta_k$ is $\mathcal{F}_k$-measurable. Hence, $V_k := \Delta_k - E[\Delta_k | \mathcal{F}_{k+1}]$ are reversed martingale differences. Since $-1 \leq \Delta_k \leq 1$, we can use Höffding–Azuma's inequality to obtain

(19) $$P\left(\sum_{k=i+1}^{j} \Delta^k - \sum_{k=i+1}^{j} E[\Delta^k | \mathcal{F}_{k+1}] < -c\right) \leq \exp\left[-\frac{2c^2}{4(j-i)}\right].$$

The inequality (18) means

$$E[\Delta_k | \mathcal{F}_{k+1}] \geq \alpha_1 - \alpha_2$$

implying that

(20) $$\sum_{k=i+1}^{j} E[\Delta^k | \mathcal{F}_{k+1}] \geq (\alpha_1 - \alpha_2)(j-i).$$

With $c = (\frac{\alpha_1 - \alpha_2}{2})(j-i)$, (19) and (20) yield

$$P\left(\sum_{k=i+1}^{j} \Delta^k < \left(\frac{\alpha_1 - \alpha_2}{2}\right)(j-i)\right)$$

$$\leq P\left(\sum_{k=i+1}^{j} \Delta^k - \sum_{k=i+1}^{j} E[\Delta^k | \mathcal{F}_{k+1}] < -\left(\frac{\alpha_1 - \alpha_2}{2}\right)(j-i)\right)$$

$$\leq \exp[-\alpha(j-i)],$$

where $\alpha = \frac{1}{2}(\frac{\alpha_1 - \alpha_2}{2})^2$. So

$$(21) \qquad P\left(\sum_{k=i+1}^{j} \Delta^k < \alpha_3(j-i)\right) \leq \exp[-\alpha(j-i)].$$

Let $E_{\Delta \text{ slope}}^n$ be the event that $\forall i, j \in I$, such that $2\varepsilon n < i < j \leq 2\varepsilon n + \sqrt{n}$ and $i + n^{0.1} \leq j$, we have

$$(22) \qquad \sum_{k=i}^{j} \Delta^k \geq \alpha_3 |i - j|.$$

By (21), for $n$ large enough, there exists a constant $c_2 > 0$ such that

$$P((E_{\Delta \text{ slope}}^n)^c) \leq n \exp[-(\alpha)n^{0.1}] \leq \exp[-c_2 \cdot n^{0.1}],$$

and hence

$$(23) \qquad P(E_{\Delta \text{ slope}}^n) \geq 1 - e^{-c_2 \cdot n^{0.1}}.$$

When the event $A_n^{\text{all}}$ holds, $E_{\text{slope}}^n$ and $E_{\Delta \text{ slope}}^n$ are equivalent. Hence

$$A_n^{\text{all}} \cap E_{\Delta \text{ slope}}^n \subset E_{\text{slope}}^n,$$

which implies

$$(24) \qquad P(E_{\text{slope}}^{nc}) \leq P((A_n^{\text{all}})^c) + P(E_{\Delta \text{ slope}}^{nc}).$$

Note that

$$(25) \quad P((A_n^{\text{all}})^c) \leq \sum_{k \in I} P(A_n^{kc}) = \sum_{k \in I} P(A_n^c | N_1 = k) \leq \sum_{k \in I} \frac{P(A_n^c)}{P(N_1 = k)},$$

where

$$(26) \qquad A_n := \{(X, Y) \in B_n\}.$$

By the local central limit theorem, there exists $c_3 > 0$ such that for all $k \in I$

$$P(N_1 = k) \geq \frac{1/c_3}{\sqrt{n}}.$$

Applying the last inequality to (25), yields

(27) $$P((A_n^{\text{all}})^c) \leq \sqrt{2}nc_3 P(A_n^c).$$

Now the inequalities (23), (27) and (24) yield

(28) $$P(E_{\text{slope}}^{nc}) \leq \sqrt{2}nc_3 P(A_n^c) + e^{-c_2 \cdot n^{0.1}}.$$

By Theorem 2.2, we have that $P(A_n^c) \leq Ce^{-c_1 n}$. Applying this to (28) gives

$$P(E_{\text{slope}}^{nc}) \leq c_3\sqrt{2}ne^{-c_1 n} + e^{-c_2 \cdot n^{0.1}},$$

which finishes the proof. □

When $E_{\text{slope}}^n$ holds, then the map

$$L : I \to \mathbb{N}$$

satisfies the conditions of Lemma 3.1 with $m = n^{0.1}$. Hence, when $E_{\text{slope}}^n$ holds, then

$$\text{Var}[L(\tilde{N}_1)] \geq \alpha_3^2 \left(1 - \frac{2n^{0.1}}{\alpha_3 \sqrt{\text{Var}[\tilde{N}_1]}}\right) \text{Var}[\tilde{N}_1].$$

Conditioning on $E_{\text{slope}}^n$, using the fact that the variance is nonnegative and $\tilde{N}_1$ and $L$ are independent,

$$E[\text{Var}[L(\tilde{N}_1)]|L(\cdot)] \geq E[\text{Var}[L(\tilde{N}_1)|E_{\text{slope}}^n]]P(E_{\text{slope}}^n)$$

$$\geq \alpha_3^2 \left(1 - \frac{2n^{0.1}}{\alpha_3 \sqrt{\text{Var}[\tilde{N}_1]}}\right) \text{Var}[\tilde{N}_1]P(E_{\text{slope}}^n).$$

Plugging the last inequality into (11) yields

(29)
$$E[\text{Var}[L(N_1)|L(\cdot)]]$$
$$\geq \alpha_3^2 \left(1 - \frac{2n^{0.1}}{\alpha_3 \sqrt{\text{Var}[\tilde{N}_1]}}\right) \text{Var}[\tilde{N}_1]P(E_{\text{slope}}^n)P(N_1 \in I).$$

By the central limit theorem, $P(N_1 \in I)$ converges to

$$P(\mathcal{N}(0,1) \in [-1,1]) > 0$$

as $n \to \infty$. [Here $\mathcal{N}(0,1)$ designate the standard normal variable.]

Note that $N_1$ is a binomial variable with parameters $2n$ and $\varepsilon$. Hence, by the central limit theorem,

$$\frac{\text{Var}[\tilde{N}_1]}{n} = \frac{\text{Var}[N_1|N_1 \in I]}{n}$$

$$\to 2\varepsilon(1-\varepsilon)P(\mathcal{N}(0,1) \in [-1,1])^{-1} \int_{-1}^{1} \phi(x)x^2\,dx,$$

where $\phi$ is the standard normal density. Together with Lemma 3.2, this implies that the right-hand side of inequality (29) divided by $n$ converges to

$$\alpha_3^2 2\varepsilon(1-\varepsilon) \int_{-1}^{1} \phi(x)x^2 \, dx > 0.$$

The inequality (9) now finishes the proof.

**4. Aligning the ones.** The rest of the paper is devoted to the proof of Theorem 2.2. The key ingredient for the following is the notation to describe the alignments. Throughout this paper we only consider alignments which align a symbol with a gap or with the same symbol in the other text. We exclude alignments which align different symbols with each other. We start with a simple example.

EXAMPLE. Take the two texts $X = 1000001$ and $Y = 1001$. The LCS of $X$ and $Y$ is $Z = 1001$. It is obtained by aligning the first one in both text and the last one and for the rest aligning as many zeros as possible. Text $X$ contains 5 zeros and text $Y$ contains 2. The maximum number of aligned zeros is thus $\min\{2,5\} = 2$. There are many alignments corresponding to the LCS $Z = 1001$. Let us present two alignments corresponding to this LCS:

| $X$ | 1 | 0 | 0 | 0 | 0 | 0 | 1 |
|---|---|---|---|---|---|---|---|
| $Y$ | 1 | 0 | 0 |   |   |   | 1 |

or another possibility:

| $X$ | 1 | 0 | 0 | 0 | 0 | 0 | 1 |
|---|---|---|---|---|---|---|---|
| $Y$ | 1 |   |   |   | 0 | 0 | 1 |

How the zeros are aligned between the ones is not important as long as we align the maximum number of zeros between the ones. Hence in general we will only describe which ones are aligned and assume that between pairs of aligned ones we align the maximum number of zeros. Let us give a further example to illustrate this. Take the sequences:

$$X = 101010101,$$
$$Y = 11010001.$$

A LCS of $X$ and $Y$ is 1101001. This LCS can be obtained with the following alignment:

(30)
| $X$ | 1 | 0 | 1 | 0 | 1 | 0 |   | 1 | 0 | 1 |
|---|---|---|---|---|---|---|---|---|---|---|
| $Y$ | 1 |   | 1 | 0 | 1 | 0 | 0 |   | 0 | 1 |

We call the portions between pairs of aligned ones *cell*.

The first cell of alignment (30) is

$$\begin{array}{c} \overline{\phantom{0}1\phantom{0}} \\ \underline{\phantom{0}1\phantom{0}} \end{array}.$$

The first cell is an exception. It is the only cell which is not comprised between two pairs of aligned ones. Instead it consists of the first pair of aligned ones and everything to its left. We only introduce this special cell in order to simplify notations later on.

The second cell of alignment (30) is

$$\begin{array}{cc} \overline{0 \quad 1} \\ 1 \end{array}.$$

The third cell of alignment (30) is

$$\begin{array}{cc} \overline{0 \quad 1} \\ 0 \quad 1 \end{array}.$$

The fourth cell of alignment (30) is

$$\begin{array}{cccc} \overline{0} & 1 & 0 & 1 \\ 0 \quad 0 & & 0 & 1 \end{array}.$$

Note that the second cell has one more zero in the $X$-part than in the $Y$-part. The third cell has the same amount of zeros in both parts. The fourth cell has two zeros in the $X$-part and three zeros in the $Y$-part. Hence the $X$-part has one zero less. The difference of zeros between the $X$-part and the $Y$-part for cell 2, 3 and 4 in this order is 1, 0 and $-1$. Cell number 1 has no zeros. Hence the difference of zeros for cell number 1 is equal to zero. Let $v_i$ denote the difference of zeros of cell $i$. We will represent alignments as the sequence of differences of zeros of their cells. For the alignment (30), this gives the representation $(v_1, v_2, v_3, v_4) = (0, 1, 0, -1)$. This sequence uniquely defines the alignment of the ones.

Let $X = X_1 \cdots X_n$ and $Y = Y_1 \cdots Y_n$ be given. As explained above, to every optimal alignment corresponds a vector $v := (v_1, \ldots, v_k)$ that shows the number of cells in the alignment $(k)$ and the difference of zeros in the cells. In every cell, the maximum amount of zeros is aligned. On the other hand, to every vector $v = (v_1, \ldots, v_k) \in \mathbb{Z}^k$ corresponds a (possible empty) family of alignments. All of them have the same pairs of aligned ones and between consecutive pairs of aligned ones, the maximum number of zeros is aligned. The alignments corresponding to $v$ can differ only in the way the zeros between aligned ones (inside cells) are aligned. Since all the alignments associated with $v$ have the same score (the same number of aligned zeros and ones), we do not care how the zeros inside a cell are aligned (as long as the maximal number of them is aligned). Therefore, in a slight imprecision we will speak of one alignment for the whole family associated with $v$. In other words, we identify each vector $v$ with an alignment. In this alignment, the number of aligned ones (cells) is $k$, the difference in the number of zero's in cell number $i$ is $v_i$ and inside a cell, the maximal number of zeros is aligned. So, in a sense, it is the "smallest" alignment which aligns exactly $k$ pairs of ones with each other and has the difference of zeros in cell $i$ equal to $v_i$, for all $i \in \{1, 2, \ldots, k\}$.

We write $|v|$ for the length of $v$. If $v \in \mathbb{R}^k$, then $|v| = k$. Let us next define rigorously the alignment associated with $v = (v_1, \ldots, v_k) \in \mathbb{Z}^k$.

DEFINITION 4.1.  Let $k \in \mathbb{N}$ and let $v = (v_1, \ldots, v_k) \in \mathbb{Z}^k$. Define $\pi(i), \nu(i)$ by induction on $i$:

- start with $\pi(0) = \nu(0) = 0$;
- for $i < k$, once $\pi(i), \nu(i)$ is defined, let $(\pi(i+1), \nu(i+1))$ be the smallest $(s, t)$ such that all of the following three conditions are satisfied.

  1. $\pi(i) < s$ and $\nu(i) < t$;
  2. $X_s = Y_t = 1$;
  3. the difference between the number of zeros of $X$ in the interval $[\pi(i), s]$ and the number of zeros of $Y$ in the interval $[\nu(i), t]$ is equal to $v_{i+1}$. Hence,

  $$v_{i+1} := \left( (s - \pi(i)) - \sum_{j=\pi(i)}^{s} X_j \right) - \left( (t - \nu(i)) - \sum_{j=\nu(i)}^{s} Y_j \right).$$

If no such $(s, t)$ exists, then $\pi(i+1) = \cdots = \pi(k) := \infty$ and $\nu(i+1) = \cdots = \nu(k) := \infty$.

The *cell* number $i$ is equal to the pair of strings:

$$C(i) := ((X_{\pi(i-1)+1}, \ldots, X_{\pi(i)}), (Y_{\nu(i-1)+1}, \ldots, Y_{\nu(i)})).$$

We define the alignment $v$ as any alignment such that the following conditions hold (provided that there exists at least one):

- $X_{\pi(i)}$ is aligned with $Y_{\nu(i)}$ for every $i = 1, \ldots, k$;
- the *number of aligned zeros* in the cell $C(i)$, denoted by $S_v(i)$, is the minimum between the number of zeros in the string $X_{\pi(i-1)+1} X_{\pi(i-1)+2} \cdots X_{\pi(i)}$ and the number of zeros in the string $Y_{\nu(i-1)+1} Y_{\nu(i-1)+2} \cdots Y_{\nu(i)}$;
- after aligning $X_{\pi(k)}$ with $Y_{\nu(k)}$, we align as many zeros as possible. Let that number be $r$.

Hence, the number of aligned zeros up to the last pair of aligned ones equal to

$$S(i) := \min\left\{ (\pi(i) - \pi(i-1)) - \sum_{j=\pi(i-1)+1}^{\pi(i)} X_j, \right.$$
$$\left. (\nu(i) - \nu(i-1)) - \sum_{j=\nu(i-1)+1}^{\nu(i)} Y_s \right\}.$$

To show that all $\pi(i), \nu(i), C(i), S(i)$ depend on $v$, we write also

$$\pi_v(i) := \pi(i), \qquad \nu_v(i) := \nu(i), \qquad C_v(i) := C(i), \qquad S_v(i) := S(i),$$
$$r_v := r.$$

We call a cell $C_v(i)$ a *u-cell*, if $v_i = u$. Thus, in a 0-cell the number of zeros in $X$-part equals the number of zeros in $Y$-part, and all the zeros in the cell are aligned. Similarly, in a 2-cell, there are 2 more zeros in the $X$-part, and these 2 zeros remain unaligned.

To summarize: every $v \in \mathbb{Z}^k$ defines an alignment. This alignment corresponds to aligning $X_{\pi_v(i)}$ with $Y_{\nu_v(i)}$, for each $i = 1, 2, \ldots, k$. These are the aligned pairs of ones: $X_{\pi_v(i)} = Y_{\nu_v(i)} = 1$. Between the aligned pairs of ones we assume that we align as many zeros as possible. Hence in cell number $i$, we align $S_v(i)$ zeros (maximum possible amount). After last pair of aligned ones, we align as many zeros as possible. The length of the common subsequence defined by alignment $v$ can now be computed as follows:

Each cell gives one aligned pair of ones. Hence, this part contributes $|v|$. Then we add for each cell the number of zeros aligned. This sums up to $\sum_{i=1}^{|v|} S_v(i)$. Finally we need to add the remaining amount of zeros $r_v$ which can be aligned but which come after the last cell. When $v \in \mathbb{Z}^k$ is such that $\pi_v(k), \nu_v(k) \leq n$, then $r_v$ is the minimum between the number of zeros in the string $X_{\pi_v(k)} \cdots X_n$ and the number of zeros in the string $Y_{\nu_v(k)} \cdots Y_n$. The length of the common subsequence defined by the alignment $v$ is now equal to

$$S_v := |v| + \sum_{i=1}^{|v|} S_v(i) + r_v.$$

The number $S_v$ is also called the *score* of the alignment $v$. This is the length of the common subsequence corresponding to $v$.

Of course, it can be that given $X = X_1 \cdots X_n$ and $Y = Y_1 \cdots Y_n$ there might not be any alignment corresponding to $v$. In this case $\pi(k) = \nu(k) = \infty$. On the other hand, if an alignment corresponding to $v$ exists, then $\pi_v(k) \leq n$ and $\nu_v(k) \leq n$. A vector $v \in \mathbb{Z}^k$ satisfying the previous condition is called *admissible*. Let $V$ designate the set of all admissible alignments, that is,

$$(31) \qquad V := \left\{ v \in \bigcup_{k>0} \mathbb{Z}^k : \pi(|v|), \nu(|v|) \leq n \right\}.$$

The set $V$, obviously, depends on $X$ and $Y$. The next statement trivially holds.

PROPOSITION 4.1.

$$(32) \qquad L_n = \max_{v \in V} \left( |v| + \sum_{i=1}^{|v|} S_v(i) + r_v \right).$$

*We say an admissible alignment $v$ is optimal if $S_v = L_n$.*

Let $v \in \bigcup_{k>0} \mathbb{Z}^k$ be nonrandom and define $|v|$ random cells $C_v(1), \ldots, C_v(|v|)$ as in Definition 4.1. One of the main advantages of defining alignments the way described above is that the cells $C_v(1), C_v(2), \ldots, C_v(|v|)$ are independent so that we can use large deviation techniques. If $v_i = v_j = u$, then, in addition to being independent, the cells $C_v(i)$ and $C_v(j)$ are both identically distributed $u$-cells. In Section 6.1, we show how to efficiently construct a $u$-cell.

## 5. The effect of changing a one into a zero.

5.1. *The events $B_n$ and $A_n$.* Recall the main idea behind Theorem 2.2: typically, when changing a randomly picked one into a zero, the score $L_n$ is likelier to increase than to decrease. More precisely, we want the conditional probability of an increase in score to be above $\alpha_1$, while the conditional probability of a decrease should be below $\alpha_2$. The constants $\alpha_1$ and $\alpha_2$, do not depend on $n$ and satisfy $\alpha_1 > \alpha_2$. By "conditional," we mean conditional on $X$ and $Y$.

EXAMPLE. Take the two texts $X = 0001000001$ and $Y = 1000010101$. An optimal alignment is given by

| $X$ | | 0 | 0 | 0 | 1 | 0 | | 0 | 0 | 0 | 0 | 1 |
|---|---|---|---|---|---|---|---|---|---|---|---|---|
| $Y$ | 1 | 0 | 0 | 0 | 0 | 1 | 0 | 1 | 0 | | | 1 |

.

The first cell in this alignment is

| | | 0 | 0 | 0 | 1 |
|---|---|---|---|---|---|
| 1 | 0 | 0 | 0 | 0 | 1 |

while the second cell is

| 0 | | 0 | 0 | 0 | 0 | 1 |
|---|---|---|---|---|---|---|
| 0 | 1 | 0 | | | | 1 |

.

Assume that the one which we switch into a zero is $Y_8$. This is a "nonaligned" one contained in the $Y$-part of cell number two. By switching $Y_8$ into a zero the LCS increases by one unit. The reason is that in cell number two, we can now align three zeros instead of only two. The new cell number two (after switching $Y_8$) looks as follows:

| 0 | 0 | 0 | 0 | 0 | 1 |
|---|---|---|---|---|---|
| 0 | 0 | 0 | | | 1 |

.

The score gets increased because $Y_8$ is on the side of the cell with strictly less zeros. We say that $Y_8$ is *on the side of a cell with less zeros.* Let us imagine next that instead of $Y_8$ the one chosen would be $X_{10}$. This one is "used" in the alignment and hence switching it could result (and does in this case) in decreasing the optimal score $L_n$ by one unit. (This is not always necessary though, as can be seen with $X_4$. When we flip $X_4$ into a zero, the score remains the same.) We call the ones which are "used" in the alignment, ones that are *matched by the alignment.* In our example, $X_4$ is matched with $Y_6$ and $X_{10}$ is matched with $Y_{10}$, $Y_8$ is not matched, nor is $Y_1$.

In the present situation, we have six ones. Each one has a probability to get picked of $1/6$. Only $Y_1$ and $Y_8$ increase the score when picked. (Here $Y_1$ is a nonmatched one on a side with more zeros. In general, such a one must not increase the score when changed into a zero. It does in this example by completely modifying the alignment and changing the number of cells.) Hence the probability of an increase in score is equal to $2/6$. Four ones, $X_4$, $X_{10}$, $Y_6$ and $Y_{10}$ could potentially decrease the score. In our example only $X_{10}$ actually does, so the conditional probability of a decrease is $1/6$. Since, in general, with longer sequences we cannot look in detail at every one, we will use as upper-bound for the probability of a decrease: the proportions of matched ones to total number of ones. In our case, this gives $4/6$ as upper bound for the probability of a decrease in score. As lower bound for the probability of an increase, we take the proportion of unmatched ones on sides with less zeros to the total number of ones. In our example, this proportion is equal to $2/6$.

From our example, it becomes clear what we need to do. We need to prove that typically there exists an optimal alignment $v$ for which:

(1) The proportion of ones that are on a side of a cell with less zeros among all ones in $X$ and $Y$ is above $\alpha_1$.

(2) The proportion of ones that are matched among all ones in $X$ and $Y$, is below $\alpha_2$.

In other words, we need to show that there exists an optimal alignment, with much less aligned ones than ones that are on a side of a cell with less zeros.

Let $N_v^-(i)$ denote the number of ones on the side with less zeros in cell number $i$. Formally, let $k \in \mathbb{N}$ and let $v = (v_1, \ldots, v_k) \in \mathbb{Z}^k$ be admissible. For $i \in [0, k]$, we define

$$N_v^-(i) := \begin{cases} 0, & \text{if } v_i = 0 \text{ (there is no side with less zeros)}; \\ \sum_{j=\nu(i)+1}^{\nu(i+1)-1} Y_j, & \text{if } v_i > 0 \text{ ($Y$ part has less zeros)}; \\ \sum_{j=\pi(i)+1}^{\pi(i+1)-1} X_j, & \text{if } v_i < 0 \text{ ($X$ part has less zeros)}. \end{cases}$$

The *total number of ones on sides with less zeros* is

$$N_v^- := \sum_{i=1}^{|v|} N_v^-(i).$$

It is important to note that $N_v^-(i)$ counts the ones inside the cell, that is, the aligned 1 that ends every cell is not counted. This means that $N_v^-(i)$ can also be zero. [In the example above, $N_v^-(1) = 0$, $N_v^-(2) = 1$ and $N^- = 1$.] Such a definition ensures that $N_v^- \geq \alpha_1 N_1$ guarantees $P(\tilde{L} - L = 1 | X, Y) \geq \alpha_1$.

Fix some constants $\alpha_1, \alpha_2$. Let $A_n$ be the event that there exists an optimal alignment $v$ such that

1. The proportions of ones on sides with less zeros is above $\alpha_1$. Hence, $N_v^- \geq \alpha_1 N_1$.
2. The proportion of aligned ones is below $\alpha_2$: $2|v| \leq \alpha_2 N_1$, where $N_1$ is the total number of ones in $X$ and $Y$.

Obviously, $A_n$ depends on the values of $\alpha_1$ and $\alpha_2$. From what we explained it follows directly that on $A_n$, the desired inequalities hold:

$$P(\tilde{L} - L = 1 | X, Y) \geq \alpha_1 \quad \text{and} \quad P(\tilde{L} - L = -1 | X, Y) \leq \alpha_2.$$

What is left to prove is that there exists $\alpha_1 > \alpha_2 > 0$ such that the event $A_n$ has probability close to one:

(33) $$P(A_n) \geq 1 - \exp[-c_1 n], \quad \text{where } c_1 > 0.$$

To be consistent with the notation in Theorem 2.2, let $B_n$ designate the set of pairs of strings $(x, y)$ for which $A_n$ holds. Hence, $(x, y) \in B_n$ if and only if

$$\{X = x, Y = y\} \in A_n.$$

We have $A_n := \{(X, Y) \in B_n\}$ and for $(x, y) \in B_n$ inequalities (4) and (5) hold.

5.2. *Breaking cells.* In the previous section we argued that we need an optimal alignment with enough ones in cell-sides with less zeros (0-cells). The problem is that many optimal alignments can have most cells with the same number of zeros on both sides. For such alignments there will also be few ones on cell-sides with less zeros. This problem is circumvented by taking an optimal alignment with most cells having same number of zeros on both sides and applying some surgery, so as to create enough cells with different numbers of zeros on the sides. This is done in such a manner that the "patient" after operation is still an optimal alignment. Let us first look at an example.

EXAMPLE. Take the texts $X = 01001001001001$ and $Y = 01010000010101$. Take the following optimal alignment

| X | 0 | 1 | 0 |   | 0 | 1 | 0 | 0 | 1 | 0 | 0 | 1 | 0 |   | 0 | 1 |
|---|---|---|---|---|---|---|---|---|---|---|---|---|---|---|---|---|
| Y | 0 | 1 | 0 | 1 | 0 |   | 0 | 0 |   | 0 | 0 | 1 | 0 | 1 | 0 | 1 |

The first cell is

| 0 | 1 |
|---|---|
| 0 | 1 |

The second cell is

| 0 |   |   | 0 | 1 | 0 | 0 | 1 | 0 | 0 | 1 |
|---|---|---|---|---|---|---|---|---|---|---|
| 0 | 1 | 0 |   | 0 | 0 |   |   | 0 | 0 | 1 |

The third cell is

| 0 |   | 0 | 1 |
|---|---|---|---|
| 0 | 1 | 0 | 1 |

All cells in the above alignment have the same number of zeros. Thus, there are no sides with less zeros and $N_v^- = 0$. Now there is a way to remedy this problem. Take cell number two. There are two ones which are "quasi" aligned: $X_5$ and $Y_4$. These two ones are only one position away from being aligned. So, if we align them, instead of the pair of zeros $X_4$ and $Y_5$, the score remains the same. When we align the pair of ones $X_5$ and $Y_4$, we split cell number two into two cells. This is how cell number two looks after this transformation:

| 0 | 0 | 1 |   | 0 | 0 | 1 | 0 | 0 | 1 |
|---|---|---|---|---|---|---|---|---|---|
| 0 |   | 1 | 0 | 0 | 0 |   | 0 | 0 | 1 |

Instead of the old cell number two, we observe the new cell number two followed by the new cell number three. The old cell number three does not change but is renamed and becomes cell number 4. The new cell number two is equal to

$$\begin{array}{|ccc|} \hline 0 & 0 & 1 \\ 0 & & 1 \\ \hline \end{array}.$$

The new cell number three is

$$\begin{array}{|cccccc|} \hline 0 & 0 & 1 & 0 & 0 & 1 \\ 0 & 0 & 0 & & 0 & 0 & 1 \\ \hline \end{array}.$$

The advantage of breaking up a cell is that the new cells have different number of zeros on each side. Hence, $N^-$ tends to increases in the process while the score remains the same. In our example, after breaking the cells, the number of ones on sides with less zeros is 1, since the new cell number three has a one on a side of less zeros. Changing this one into a zero will increase the score. The breaking up process helps up get rid of the problem of having too many cells with the same number of zeros on both sides. Note that the breaking up the cell does not necessary increase the number $N^-$: although after breaking a cell, both new cells have different number of zeros, it might happen that both of them have no ones on the side of less zeros. In this case, the number $N^-$ does not increase. However, once we have an optimal alignment with enough nonzero cells, the probability is high to also find enough ones on sides with less zeros.

Let us define what we saw in the previous numerical example in a precise fashion.

DEFINITION 5.1. Let $k \in \mathbb{N}$, $v \in \mathbb{Z}^k \cap V$, $i \leq k$ and $v_i = 0$. We say that cell $i$ of $v$ can be *broken up* if there exists $j$ and $j'$ satisfying all of the following:

1. $X_j = Y_{j'} = 1$.
2. $\pi(i) < j < \pi(i+1)$ and $\nu(i) < j' < \nu(i+1)$.
3. The difference between the number of zeros in the strings

$$X_{\pi(i)+1} X_{\pi(i)+2} \cdots X_{j-1} \quad \text{and} \quad Y_{\nu(i)+1} Y_{\nu(i)+2} \cdots Y_{j'-1}$$

is one or minus one. Hence

$$1 = \left| \left( j - \pi(i) - \sum_{l=\pi(i)+1}^{j} X_l \right) - \left( j' - \nu(i) - \sum_{l=\nu(i)+1}^{j'} Y_l \right) \right|.$$

5.3. *Optimal alignment contained in $V_n$.* Recall that $Y_i$ and $X_i$ are i.i.d. Bernoulli random variables with parameter $\varepsilon$. In Section 6, we will show that with high probability $L_n$ is larger by $0.1\varepsilon^2 n$ than half of the total amount of zeros in $X$ or in $Y$. Let us briefly explain the use of this fact. When

$$L_n \geq \frac{N_0}{2} + a,$$

where $N_0$ is the total number of zeros in $X$ and $Y$ and $a > 0$, there are two immediate consequences:

(1) In any optimal alignment $v$ there need to be at least $a$ pairs of aligned ones. Hence, any optimal alignment $v$ needs to be contained in the set $\bigcup_{k \geq a} \mathbb{Z}^k$.

(2) Any optimal alignment $v$ in $\mathbb{Z}^k$, satisfies

$$\sum_{i=1}^{k} |v_i| \leq 2k. \tag{34}$$

Otherwise the unmatched zeros (at least $\sum_{i=1}^{k} |v_i|$) would out-number the aligned ones (the number of aligned ones is $2k$) bringing the score below an alignment with only zeros aligned. Indeed, the number of nonaligned zeros in the alignment $v$ is at least $\sum_{i=1}^{k} |v_i|$, so the number of aligned zeros is at most

$$\frac{N_0 - \sum_{i=1}^{k} |v_i|}{2}$$

and (34) follows from the inequalities

$$\frac{N_0}{2} < L_n \leq \frac{N_0 - \sum_{i=1}^{k} |v_i|}{2} + k.$$

When we take $0.1\varepsilon^2 n$ for $a$, conditions (1) and (2) can be expressed by saying that any optimal alignment $v$ is necessarily contained in the set $V_n$, where

$$V_n := \bigcup_{k \geq 0.1\varepsilon^2 n} V(k), \tag{35}$$

and $V(k) \subset \mathbb{Z}^k$ is defined as follows

$$V(k) := \{(v_1, v_2, \ldots, v_k) \in \mathbb{Z} | |v_1| + \cdots + |v_k| \leq 2k\}. \tag{36}$$

The fact that any optimal alignment is typically contained in $V_n$ is very useful. The set $V_n$ is relatively small [see the bound (50)]. So, whenever we want to prove the likeliness of a property for the optimal alignment, we prove the property to hold typically for every alignment in $V_n$. The tremendous advantage of this approach is that for every (nonrandom) $v \in V_n$, the alignment associated with $v$ has a simple distribution: the cells are independent. This allows us to use large deviation techniques. In contrast, in the optimal alignment the cells are correlated in a complex and poorly understood manner.

A cell which has different number of zeros in its $X$-part and in its $Y$-part is called a *nonzero cell*. We say that an alignment $v \in \mathbb{Z}^k$ has more than 1% nonzero cells if

$$|\{i \in [1,k] | v_i \neq 0\}| \geq 0.01k.$$

Let $V_{1\%}$ be the subset of $V_n$ consisting of the alignments which have at least 1% of nonzero cells, that is,

$$V_{1\%} := \{v \in V_n | v \text{ has more than } 1\% \text{ nonzero cells}\}.$$

Let

$$V_{1\%}^c := V_n - V_{1\%}.$$

5.4. *The events.* Recall that for a vector $v$ we associate $|v|$ random cells $C_v(1), \ldots, C_v(|v|)$ defined as a function of random i.i.d., Bernoulli random sequences $X_1, X_2, \ldots$ and $Y_1, Y_2, \ldots$. In the following we define some events that capture the typical behavior of these random cells.

Recall that $N_1$ denotes the total number of ones in $X$ and $Y$, $N_1 = \sum_{i=1}^n (X_i + Y_i)$. Let $v$ be an admissible alignment, that is, $v \in V$ or, equivalently, $\pi_v(|v|), \nu_v(|v|) \leq n$.

Let $N_{1v}$ designate the number of ones up to the last cell of $v$:

$$N_{1v} := \left( \sum_{j=1}^{\pi_v(|v|)} X_j + \sum_{j=1}^{\nu_v(|v|)} Y_j \right).$$

Finally, we define the number of ones after the last cell

$$R_v = \sum_{j=\pi(|v|)+1}^n X_j + \sum_{j=\nu(|v|)+1}^n Y_j.$$

DEFINITION 5.2.
- Let $E_4$ designate the event that every optimal alignment belongs to the set $V_n$.
- Let $D$ be the event that for all $v \in V_{1\%}^c$, at least 1% of the cells can be broken up. So,

$$D := \bigcap_{v \in V_{1\%}^c} D_v,$$

where $D_v$ is the event that at least 1% of the cells $C_v(1), \ldots, C_v(|v|)$ can be broken up.

- Let $F$ be the event that every $v \in V_{1\%}$ has at least $2\alpha_1\%$ of ones in $C_v(1), \ldots, C_v(|v|)$ on a side of less zeros. Hence,

$$F := \bigcap_{v \in V_{1\%}} F_v,$$

where $F_v$ is the event that

$$N_v^- \geq 2\alpha_1 N_{1v}.$$

- Let $G$ be the event that every $v \in V_{1\%}$ has no more than $\alpha_2\%$ of matched ones. Hence

$$G := \bigcap_{v \in V_{1\%}} G_v,$$

where $G_v$ is the event that

$$2|v| \leq \alpha_2 N_{1v}.$$

- Let $K$ be the event that there exists an optimal alignment $v$ such that

$$R_v \leq N_{1v}.$$

In the next section, we shall prove that all the defined events hold with high probability. Note the importance of the breaking up notion. The events $F$ and $G$ together with the event $K$ basically prove (4) and (5) for the case when the optimal alignment has at least 1% nonzero cells, that is, it belongs to $V_{1\%}$. But every optimal alignment needs not belong to $V_{1\%}$. However, the event $D$ ensures that for every alignment from $V_{1\%}^c$, there exists another alignment $v' \in V_{1\%}$ *with the same score*. So, when the events $E_4$ and $D$ both hold, then there exists an optimal alignment in $V_{1\%}$. To this optimal alignment we can apply $F$, $G$ and $K$ and get the inequalities (4) and (5). These considerations lead to the next lemma, which is our main combinatorial lemma. Recall the definition of $A_n$ in Section 5.1.

LEMMA 5.1.

(37) $$E_4 \cap D \cap F \cap G \cap K \subset A_n.$$

PROOF. Recall that $A_n$ holds if there exists an optimal alignment, say $w$, such that the following conditions are fulfilled:

(1) the proportion of ones with less zeros is above $\alpha_1$: $N_w^- \geq \alpha_1 N_1$;
(2) the proportion of aligned ones is below $\alpha_2$: $2|w| \leq \alpha_2 N_1$.

By the event $K$ we know that there exists an optimal alignment $v$ such that $R_v \leq N_{1v}$. When $E_4$ holds, then $v$ is contained in the set $V_n$. Assume that $v$ contains less than 1% of cells with different number of zeros on their sides, that is, $v \in V^c_{1\%}$. Then, the event $D$ ensures that we can break up $v$ so that it gets more than 1% of nonzero cells and still remains optimal. Let that alignment be $w$. By doing the break up, the number of ones after the last cell remains unchanged, that is, $R_w = R_v$. Moreover, breaking up only increases the number of aligned ones, so $N_{1v} \leq N_{1w}$. Hence, there exists an optimal $w \in V_{1\%}$ such that $R_w \leq N_{1w}$. The events $F$ applies to $w$. Hence

$$(38) \qquad \frac{N^-_w}{N_{1w}} = \frac{N^-_w}{N_1} \cdot \frac{N_1}{N_{1w}} \geq 2\alpha_1.$$

Since $w$ is admissible, $N_1 = N_{1w} + R_w$. Hence

$$\frac{N_1}{N_{1w}} = \frac{N_{1w} + R_w}{N_{1w}} \leq 2.$$

Using the last equality with (38) yields

$$\frac{N^-_w}{N_1} \geq \alpha_1.$$

This is the first statement on the event $A_n$.

Since $w \in V_{1\%}$, the event $G$ guarantees that there is a proportion of less than $\alpha_2\%$ matched ones: $2|w| \leq \alpha_2 N_{1w} \leq \alpha_2 N_1$. This proves the second statement of the event $A_n$. □

5.5. *Proof of Theorem 2.2.* From (37) it follows that

$$(39) \qquad P(A^c_n) \leq P(E^c_4) + P(D^c) + P(F^c) + P(G^c) + P(K^c).$$

So, the proof of Theorem 2.2 is accomplished, if we show that there exists $\alpha_1 > \alpha_2 > 0$ and $\varepsilon_0 > 0$ such that the events $P(E^c_4)$, $P(D^c)$, $P(F^c)$, $P(G^c)$ and $P(K^c)$ are exponentially small in $n$, provided $\varepsilon \leq \varepsilon_0$. In Lemma 7.9, we prove the existence of constants $\alpha_1 > 0$ and $C_F$, not depending on $\varepsilon$, as well as a constants $c_F(\varepsilon)$ and $\varepsilon_1 > 0$ such that $P(F^c) \leq C_F \exp[-c_F n]$, if $\varepsilon < \varepsilon_1$. In Lemma 7.10, we prove that for every $0 < \alpha_2 < \alpha_1$, there exists $\varepsilon_2 < \varepsilon_1$, depending on $\alpha_2$, such that for every $\varepsilon \leq \varepsilon_2$, $P(G^c) \leq C_G \exp[-c_G n]$, where $C_G$ and $c_G$ are some constants (possibly depending on $\varepsilon$). In Lemmas 6.2, 7.4 and 7.11, we prove the existence of $\varepsilon_3 > 0$, finite constants $c_E, c_D, c_K$ as well as $C_E, C_D, C_K$, possibly depending on $\varepsilon$, such that

$$P(D^c) \leq C_D \exp[-c_D n], \qquad P(E^c_4) \leq C_E \exp[-c_E n],$$
$$P(K^c) \leq C_K \exp[-c_K n],$$

provided $\varepsilon < \varepsilon_3$. Thus, if $\varepsilon < \varepsilon_0 := \min\{\varepsilon_1, \varepsilon_2, \varepsilon_3\}$, all the events $P(E^c_4)$, $P(D^c)$, $P(F^c)$, $P(G^c)$, $P(K^c)$ have exponentially small probabilities.

The proofs that $D^c$, $F^c$, $G^c$ and $K^c$ all have exponentially small probability in $n$ uses the representation of alignments as elements of $V_n$. All these events state that a certain property holds for every alignment in $V_n$. The proof that they have high probability goes as follows: for a given nonrandom alignment $v \in V_n$, the cells are independent. Hence, one can use large deviation techniques. It then only remains to prove that the large deviation rate beats the number of elements in the set $V_n$.

## 6. Preliminary bounds.

6.1. *A useful approach.* In the sequel, we are often going to use the following way of constructing random sequences $X_1, X_2, \ldots$ and $Y_1, Y_2, \ldots$. Let $\xi_1, \xi_2, \ldots$ be a sequence of i.i.d. random variables with the distribution of $\xi$ being following:

$$P(\xi = 0) = 1 - \varepsilon, \qquad P(\xi = 1) = \varepsilon(1 - \varepsilon), \qquad \ldots,$$

$$P(\xi = n) = \varepsilon^n (1 - \varepsilon), \qquad \ldots.$$

The distribution of $\xi_i + 1$ is geometric. The random variables $\xi_i$ stand for the number of 1's between the 0's: $\xi_1$ is the number of ones before the first 0, $\xi_2$ is the number of ones between the first and second 0 and so on. For example, if $(\xi_1, \xi_2, \xi_3, \xi_4, \xi_5, \xi_6) = (0, 2, 0, 0, 1, 0)$, then before the first 0, there are no ones; between first and second zero, there are 2 ones; between second an third zero, there are again no ones, and so on. Hence, the corresponding sequence $X_1, X_2, \ldots$ begins with $0, 1, 1, 0, 0, 0, 1, 0, 0, \ldots$. Similarly, with the help of the random variables $\eta_1, \eta_2, \ldots$, we construct the sequence $Y_1, Y_2, \ldots$.

Recall our task: we are given a fixed vector $v = (v_1, \ldots, v_k)$ and we aim to construct the random cells (using i.i.d. random sequences $X$ and $Y$) $C_v(i)$ as in Definition 4.1: at first we wait for the first time such that a pair of ones can be aligned so that the difference of zeros between $X$- and $Y$-part is $u = v_1$ [so we get $C_v(1)$]; then we start afresh with $u = v_2$ and so on. In terms of $\xi$ and $\eta$ variables, it is relatively easy. Indeed, to get a 0-cell, we look for the smallest time $i$ such that $\xi_i \neq 0$, $\eta_i \neq 0$. So, a 0-cell can be constructed using the stopping time $T$, where

$$(40) \qquad T := \min\{i = 1, 2, \ldots : \xi_i \neq 0, \eta_i \neq 0\}.$$

To get a $-u$ cell ($u > 0$), we look for the smallest time $T$ such that $\xi_i \neq 0$ and $\eta_{u+i} \neq 0$. Hence, a $-u$-cell is constructed using the stopping time $T$, where

$$(41) \qquad T := \min\{i = 1, 2, \ldots : \xi_i \neq 0, \eta_{u+i} \neq 0\}.$$

In other words a cell with $v_i = u$ can be viewed in the following way: we first set $u$ zeros aside on side $X$ if $u \geq 0$ and on side $Y$ otherwise. Then we align consecutive pairs of zeros, until we meet for the first time a pair of aligned zeros both directly followed by a one. Let us look at a numerical example:

EXAMPLE. Take $v_1 = u = 2$. Let $X = 000101\ldots$ and $Y = 001\ldots$. We put aside the first two zeros in $X$. From there, we align all the zeros until we meet two zeros both followed directly by a one. Here, this gives the cell

| $X$ | 0 | 0 | 0 | 1 | 0 | 1 |
|---|---|---|---|---|---|---|
| $Y$ | | | 0 | | 0 | 1 |

6.2. *A bound on $L_n$.* A rough lower bound for the typical length of the LCS, is obtained as follows.

1. First only align all the zeros you can. You get approximately a common subsequence of length $(1-\varepsilon)n$ consisting only of zeros.
2. Having aligned as many zeros as you could in 1, take the ones which can be aligned without disturbing the already aligned zeros. In terms of previous subsection, it gives (approximatively) additional

$$\sum_{i=1}^{(1-\varepsilon)n} \min\{\xi_i, \eta_i\}$$

ones, since between $i-1$st and $i$th pair of aligned zeros, there are $\xi_i$ ones in $X$ and $\eta_i$ ones in $Y$. The random variables $\xi_i + 1$ and $\eta_i + 1$ are Geometrically distributed with parameter $(1-\varepsilon)$. This means that $\min\{\xi_i, \eta_i\} + 1 \sim G(1-\varepsilon^2)$, so

$$E\min\{\xi_i, \eta_i\} = \frac{1}{1-\varepsilon^2} - 1.$$

So, in average the ones contribute $(\frac{1}{1-\varepsilon^2} - 1)(1-\varepsilon)n$.

In the way described above we get a common subsequence of length about

$$(42) \quad \left[\left(\frac{1}{1-\varepsilon^2} - 1\right)(1-\varepsilon) + (1-\varepsilon)\right]n = \left[(1-\varepsilon) + \frac{\varepsilon^2}{1+\varepsilon}\right]n = \frac{n}{1+\varepsilon}.$$

To stay on the safe side, we bound $L_n$ by a quantity that is little smaller than (42); we take $[(1-\varepsilon) + 0.9\varepsilon^2]n$.

Let $E$ denote the event that the LCS is longer than $((1-\varepsilon) + 0.9\varepsilon^2)n$, that is,

$$E := \{L_n \geq ((1-\varepsilon) + 0.9\varepsilon^2)n\}.$$

LEMMA 6.1. *There exists $\varepsilon_3 > 0$ such that for every $\varepsilon < \varepsilon_3$*

$$P(E) \geq 1 - 5e^{-an},$$

*where $a(\varepsilon) > 0$.*

PROOF. Let $\delta \in (0, 0.5)$. Define the events (they depend on $\delta$)

$$E_2^x := \left\{ \left| \sum_{i=1}^n X_i - n\varepsilon \right| \leq \delta \varepsilon n \right\}, \qquad E_2^y := \left\{ \left| \sum_{i=1}^n Y_i - n\varepsilon \right| \leq \delta \varepsilon n \right\}.$$

When $E_2^x$ holds, then $X_1, \ldots, X_n$ has at least $(1 - (1+\delta)\varepsilon)n$ zeros. On $E_2^y$, the same holds for $Y_1, \ldots, Y_n$. Let

$$E_2(\delta) := E_2^x \cap E_2^y.$$

When $E_2$ holds, then at least $(1 - (1+\delta)\varepsilon)n$ zeros can be aligned.

Let $\zeta_i := \min\{\xi_i, \eta_i\}$, where $\xi_i$, $\eta_i$ are i.i.d. random variables, $\xi_i + 1 \sim G(1-\varepsilon)$. So, $\zeta_i + 1 \sim G(1-\varepsilon^2)$. From Proposition A.1 (see the Appendix), it follows that for every $\alpha < 1$,

$$(43) \qquad P\left( \sum_{i=1}^m \zeta_i < \left( \frac{\alpha}{1-\varepsilon^2} \right) m - m \right) \leq e^{-C(\alpha)m},$$

where $C(\alpha) = \alpha - 1 - \ln \alpha$. Let $E_1$ be the event that $\sum_{i=1}^m \zeta_i$ is at least $\frac{m\alpha}{1-\varepsilon^2} - m$, where $m(\delta) := (1 - (1+\delta)\varepsilon)n$ and $\alpha < 1$. So

$$E_1(\alpha, \delta) := \left\{ \sum_{i=1}^m \zeta_i \geq \left( \frac{\alpha}{1-\varepsilon^2} \right) m \right\}.$$

When $E_1$ and $E_2$ both hold, then one can align $m$ zeros and $\frac{m\alpha}{1-\varepsilon^2} - m$ ones between them. This means that on $E_1 \cap E_2$, the length of the longest common subsequence has the following lower bound

$$L_n \geq \frac{m\alpha}{1-\varepsilon^2} = \alpha \left( \frac{1 - \varepsilon(1+\delta)}{1-\varepsilon^2} \right) n = \alpha \left( \frac{1}{1+\varepsilon} - \frac{\delta\varepsilon}{1-\varepsilon^2} \right) n$$

$$= \alpha \left( 1 - \varepsilon + \frac{\varepsilon^2}{1+\varepsilon} - \frac{\delta\varepsilon}{1-\varepsilon^2} \right) n.$$

Let us compare the right-hand side of the previous inequality with $((1-\varepsilon) + 0.9\varepsilon^2)n$. Since

$$(44) \quad \begin{aligned} &\alpha \left( 1 - \varepsilon + \frac{\varepsilon^2}{1+\varepsilon} - \frac{\delta\varepsilon}{1-\varepsilon^2} \right) - (1-\varepsilon) - 0.9\varepsilon^2 \\ &= (1-\varepsilon)(\alpha - 1) + \varepsilon^2 \left( \frac{\alpha}{1+\varepsilon} - 0.9 \right) - \frac{\alpha\delta\varepsilon}{1-\varepsilon^2}, \end{aligned}$$

we see that for $\alpha = 1 - \varepsilon^3$ and $\delta = \varepsilon^2$, (44) is positive, provided $\varepsilon$ is small enough. So, if $\alpha = 1 - \varepsilon^3$ and $\delta = \varepsilon^2$, there exists $\varepsilon_3 > 0$ such that for every $\varepsilon < \varepsilon_3$, $E_1 \cap E_2 \subset E$, implying that

$$P(E^c) \leq P(E_2^c(\varepsilon^2)) + P(E_1^c(1-\varepsilon^3, \varepsilon^2)).$$

Finally, let us bound the probabilities. For any $\delta > 0$, from Höffding's inequality, it follows that

$$P((E_2^x)^c) \leq 2\exp[-2(\delta\varepsilon)^2 n], \qquad P((E_2^y)^c) \leq 2\exp[-2(\delta\varepsilon)^2 n].$$

From (43), we get that

$$P(E_1^c(\alpha,\delta)) \leq \exp[-C(\alpha)m(\delta)] = \exp[-C(\alpha)(1-(1+\delta)\varepsilon)n].$$

Take $\alpha = 1 - \varepsilon^3$ and $\delta = \varepsilon^2$ to obtain that, for every $\varepsilon < \varepsilon_3$

$$P(E^c) \leq 4\exp[-2\varepsilon^6 n] + \exp[-C(1-\varepsilon^3)(1-(1+\varepsilon^2)\varepsilon)n] \leq 5e^{-an},$$

where $a(\varepsilon) = \min\{2\varepsilon^6, C(1-\varepsilon^3)(1-(1+\varepsilon^2)\varepsilon)\}$. □

Note that Lemma 6.1 gives a lower bound for the Chvatal–Sankoff constant: $\frac{1}{1+\varepsilon}$.

If $0 < \delta \leq 0.8\varepsilon$, then on $E_2$

$$(45) \qquad N_0^x \leq n[(1-\varepsilon) + 0.8\varepsilon^2], \qquad N_0^y \leq n[(1-\varepsilon) + 0.8\varepsilon^2],$$

where $N_0^x$ and $N_0^y$ are the number of zeros in $X$ and $Y$, respectively. In this case, hence,

$$(46) \qquad \frac{N_0}{2} \leq n[(1-\varepsilon) + 0.8\varepsilon^2],$$

where $N_0$ is the number of zeros in $X$ and $Y$. On the other hand, if $E$ holds, then

$$(47) \qquad L_n \geq n[(1-\varepsilon) + 0.9\varepsilon^2].$$

So, if $0 < \delta \leq 0.8$ and $E_2 \cap E$ holds, then

$$(48) \qquad \frac{N_0}{2} + (0.1)\varepsilon^2 n \leq L_n.$$

As explained in subsection 5.3, (48) implies (34), that is, $\sum_{i=1}^k |v_i| \leq 2k$. We also showed that (48) implies that in any optimal alignment there are at least $(0.1)\varepsilon^2 n$ pairs of aligned ones. Thus, if $0 < \delta \leq 0.8$ and $E_2 \cap E$ hold, then any optimal alignment must belong to $V_n$, where the set of alignments $V_n$ has been defined in (35). Recall that $E_4$ designates the event that every optimal alignment belongs to $V_n$.

LEMMA 6.2. *There exists $\varepsilon_3 > 0$ such that for $\varepsilon < \varepsilon_3$, it holds,*

$$P(E_4) \geq 1 - 5\exp[-an] - 4\exp[-2(0.8\varepsilon)^2 \varepsilon n].$$

PROOF. We saw that $E_2(0.8\varepsilon) \cap E \subset E_4$. Proposition 4.1 now finishes the proof. □

## 7. Bounding the probabilities.

7.1. *Combinatorics.* In the following, we use $C_n^k$ as the number of combinations, that is,

$$C_n^k = \frac{n!}{k!(n-k)!}.$$

We also make use of the following fact: the number of $m$-dimensional vectors with nonnegative entries summing up to exactly $n$ is $C_{n+m-1}^{m-1}$. To see this, use the induction by $m$: for $m = 2$, it trivially holds. Suppose that the formula folds for $m$. Consider the $m+1$-dimensional vectors with nonnegative entries summing up exactly $n$. The first $m$ components determine the vectors; since the first $m$ components can sum up to $0, \ldots, n$, by the assumption the number of $m+1$-dimensional vectors with nonnegative entries summing up exactly $n$ is

(49) $$C_{m-1}^{m-1} + C_{1+m-1}^{m-1} + C_{2+m-1}^{m-1} + \cdots + C_{n+m-1}^{m-1}.$$

Using the fact that $C_{k+m-1}^{m-1} + C_{k+m-1}^{m} = C_{k+m}^{m}$, it is easy to show that (49) equals to $C_{n+m}^{m}$.

LEMMA 7.1. *For $k \geq 1$, we have*

(50) $$|V(k)| \leq 2^k C_{3k}^k \leq 16^k.$$

PROOF. Let

$$V^+(k) = \{(v_1, \ldots, v_k) \in (\mathbb{Z}^+)^k : v_1 + \cdots + v_k \leq 2k\},$$

where $\mathbb{Z}^+ = \{0, 1, \ldots\}$. Thus, $|V^+(k)|$ is the number of $k$-dimensional vectors with nonnegative integer entries and summing up to at most $2k$. By adding one more component, we get that $|V^+(k)|$ is equal to the number of $k+1$-dimensional vectors with nonnegative integer entries and summing up to exactly $2k$. The number of such vectors is $C_{2k+k+1-1}^{k+1-1} = C_{3k}^k$. It follows that

$$|V^+(k)| = C_{3k}^k \leq 2^{3k}.$$

(Here $2^{3k}$ represents the number of subsets of a set of size $3k$. Of course this upper bound is far from being optimal, but it is still sufficient for our purpose.) For every $k$-dimensional vector, there are at most $2^k$ ways to assign the signs of the entries. This then yields

$$|V(k)| \leq 2^k C_{3k}^k \leq 2^{4k} = 16^k. \qquad \square$$

Let

$$I(v_1, \ldots, v_k) = |\{i \in \{1, \ldots, k\} : v_i \neq 0\}|.$$

LEMMA 7.2.

(51) $$|V^c_{1\%}(k)| \leq \exp[0.1262k],$$

*where*

$$V^c_{1\%}(k) := V(k) \cap \{(v_1, \ldots, v_k) \in \mathbb{Z}^k : I(v_1, \ldots, v_k) \leq 0.01k\}.$$

PROOF. Without loss of generality assume that $0.01k$ is an integer. Consider the set of $0.01k$-dimensional vectors with nonnegative integer entries and summing up to at most $2k$. Let this set be

$$W^+(k) := \left\{(w_1, \ldots, w_{0.01k}) \in \mathbb{Z}^{+0.01k} : \sum_{i=1}^{0.01k} w_i \leq 2k\right\}.$$

We know that

(52) $$|W^+(k)| = C^{0.01k+1-1}_{2k+0.01k+1-1} = C^{0.01k}_{2.01k} = C^{0.01/2.01(2.01)k}_{2.01k}.$$

In order to bound the number of combinations $C^{qm}_m$, where $q \in (0,1)$ (and $qm$ is an integer), we note that

(53) $$C^{qm}_m \leq q^{-qm}(1-q)^{-m(1-q)},$$

that follows from the fact that $C^{mq}_m \cdot q^{qm}(1-q)^{m(1-q)} \leq 1$. Using (53) with $m = 2.01k$ and $q = \frac{0.01}{2.01}$, (52) yields

$$|W^+(k)| \leq \left(\frac{2.01}{0.01}\right)^{0.01k} \left(\frac{2.01}{2}\right)^{2k} = (201)^{0.01k}(1.005)^{2k}.$$

Here are $2^{0.01k}$ ways to assign the signs, so

$$|W(k)| \leq 2^{0.01k}(201)^{0.01k}(1.005)^{2k} = (402)^{0.01k}(1.005)^{2k},$$

where

$$W(k) := \left\{(w_1, \ldots, w_{0.01k}) \in \mathbb{Z}^{0.01k} : \sum_{i=1}^{0.01k} |w_i| \leq 2k\right\}.$$

Obviously,

$$|V^c_{1\%}(k)| \leq C^{0.01k}_k |W(k)|.$$

With (53), we have

$$C^{0.01k}_k \leq (100)^{0.01k}\left(\frac{100}{99}\right)^{0.99k},$$

implying that

$$|V^c_{1\%}(k)| \leq (40{,}200)^{0.01k}(1.005)^{2k}\left(\frac{100}{99}\right)^{0.99k}$$

$$< 1.1345^k < \exp[0.1262k]. \qquad \square$$

7.2. *The event D*. Recall that $D_v$ denotes the event that 1% of the cells of $v$ can be broken up. Note that the following bound holds for every $\varepsilon \in (0, \frac{1}{2}]$.

LEMMA 7.3. *Let $v \in V_{1\%}^c(k)$. Then*

(54) $$P(D_v^c) \leq \exp[-0.16k].$$

PROOF. Let us calculate the probability that a 0-cell is breakable. We use the approach introduced in Section 6.1. Recall the definition of $T$ in (40). With this construction, being breakable means the existence of $1 \leq i \leq T$ such that

$$\xi_i \neq 0, \qquad \eta_i = 0, \qquad \xi_{i+1} = 0, \qquad \eta_i \neq 0$$

or

$$\xi_i = 0, \qquad \eta_i \neq 0, \qquad \xi_{i+1} \neq 0, \qquad \eta_i = 0.$$

Let

$$U_1 := \min\{i = 2, \ldots : \xi_{i-1} \neq 0, \eta_{i-1} = 0, \xi_i = 0, \eta_i \neq 0\},$$
$$U_2 := \min\{i = 2, \ldots : \xi_{i-1} = 0, \eta_{i-1} \neq 0, \xi_i \neq 0, \eta_i = 0\},$$
$$U := U_1 \wedge U_2.$$

Let

$$\mathcal{X} := \{0, 1, 2, \ldots\}, \qquad \mathcal{X}^+ := \{1, 2, \ldots\}.$$

With those stopping times, the probability that a 0 cell is breakable is $P(U < T)$. Let us estimate it (from below). An easy way is to consider the disjoint pairs of indexes $(1, 2), (3, 4), \ldots, (2j-1, 2j), \ldots$ and restrict the stopping time $U$ to the even integers only. So, we define the independent random vectors

$$Z_j = (\xi_{2j-1}, \eta_{2j-1}, \xi_{2j}, \eta_{2j}), \qquad j = 1, 2, \ldots,$$
$$U_1' := \min\{j = 1, 2, \ldots : \xi_{2j-1} \neq 0, \eta_{2j-1} = 0, \xi_{2j} = 0, \eta_{2j} \neq 0\}$$
$$= \min\{j = 1, 2, \ldots : Z_j \in A_1\},$$
$$U_2' := \min\{i = 1, 2, \ldots : \xi_{2j-1} = 0, \eta_{2j-1} \neq 0, \xi_{2j} \neq 0, \eta_{2j} = 0\}$$
$$= \min\{j = 1, 2, \ldots : Z_j \in A_2\},$$
$$U' := U_1' \wedge U_2' = \min\{j = 1, 2, \ldots : Z_j \in A_2 \cup A_1\},$$
$$T' := \{j = 1, 2, \ldots : Z_j \in B_1 \cup B_2\},$$

where

$$A_1 := \mathcal{X}^+ \times \{0\} \times \{0\} \times \mathcal{X}^+, \qquad A_2 := \{0\} \times \mathcal{X}^+ \times \mathcal{X}^+ \times \{0\},$$
$$B_1 = \mathcal{X}^+ \times \mathcal{X}^+ \times \mathcal{X} \times \mathcal{X}, \qquad B_2 = \mathcal{X} \times \mathcal{X} \times \mathcal{X}^+ \times \mathcal{X}^+.$$

Clearly,

$$2U' - 1 \geq U, \quad 2T' - 1 \leq T,$$
$$P(U < T) \geq P(2U' - 1 < T) \geq P(2U' - 1 < 2T' - 1) = P(U' < T').$$

Since the random variables $Z_j$ are independent, the probability of the right-hand side is easy to calculate:

$$P(U' < T') = \frac{P(Z_1 \in A_2 \cup A_1)}{P(Z_1 \in A_2 \cup A_1) + P(Z_1 \in B_2 \cup B_1)} = \frac{2\varepsilon^2(1-\varepsilon)^2}{2\varepsilon^2(1-\varepsilon)^2 + 2\varepsilon^2 - \varepsilon^4}$$
$$= \frac{2(1-\varepsilon)^2}{2(1-\varepsilon)^2 + 2 - \varepsilon^2}.$$

It is easy to check that the function

$$\varepsilon \mapsto q(\varepsilon) := \frac{2(1-\varepsilon)^2}{2(1-\varepsilon)^2 + 2 - \varepsilon^2}$$

is decreasing in $[0, \frac{1}{2}]$, which implies

$$q(\varepsilon) \geq \frac{2(1/2)^2}{2(1/2)^2 + 2 - (1/2)^2} = \frac{2}{9}.$$

Let $v = (v_1, \ldots, v_k) \in V_{1\%}^c$. This means that the number of zero cells $m$ is at least $0.99k$. Let $J$ be the index set of zero cells and let for every $j \in J$, $I_j$ be the Bernoulli variable that is one if and only if the cell $v_j$ is breakable. Clearly, the random variables $I_j$ are and for every $\varepsilon > 0$, $p(\varepsilon) := P(I_j = 1) \geq q(\varepsilon) \geq \frac{2}{9}$.

In the following, we use the following result: let $Z$ be a binomial random variable with parameters $p$ and $m$. Let $0 < a < p$. Then

$$(55) \quad P(Z < am) \leq \left(\frac{p}{a}\right)^{am} \left(\frac{1-p}{1-a}\right)^{(1-a)m} \leq \left(\frac{p}{a}\right)^{am} \exp[(a-p)m]$$

(see, e.g., [8], page 130). Using (55) and the facts that $p := p(\varepsilon) \geq \frac{2}{9}$ as well as $m \geq 0.99k$, we get

$$P(D_v^c) = P\left(\sum_{j \in J} I_j < 0.01m\right) \leq (100p)^{0.01m} \exp[(0.01 - p)m]$$
$$\leq \exp[(0.01 \ln 100 + (0.01 - p))m]$$
$$\leq \exp[(0.047 + (0.01 - p))0.99k] \leq \exp[-0.16k]. \quad \square$$

LEMMA 7.4. *There exists $C_D < \infty$ such that*

$$(56) \quad P(D^c) \leq C_D \exp[-0.0438(0.1\varepsilon^2)n].$$

Proof.
$$D(k) := \bigcap_{v \in V_{1\%}^c(k)} D_v.$$

With (51) and (54), we get

$$P(D^c(k)) \leq \sum_{v \in V_{1\%}^c(k)} P(D_v^c) \leq \exp[(-1.16 + 0.1262)k] = \exp[-0.0438n].$$

Since $k \geq (0.1\varepsilon^2)n$, we find:

$$P(D^c) \leq \sum_{k \geq (0.1\varepsilon^2)n} P(D^c(k)) \leq \sum_{k \geq (0.1\varepsilon^2)n} \exp[-0.0438k]$$
$$= C_D \exp[-0.0438(0.1\varepsilon^2)n],$$

where

$$C_D := (1 - \exp[-0.0438])^{-1}. \qquad \square$$

7.3. *The event $F$.* The following large deviation result is proven in the [Appendix](#).

LEMMA 7.5 (Large deviation for geometric random variables). *Let $G_1, \ldots, G_m$ be i.i.d. random variables with geometric distribution $G(p)$. There exists $0 < \alpha_0 < 1$, not depending on $p$, such that for every $\alpha \leq \alpha_0$, the inequality*

$$(57) \qquad P\left(\sum_{i=1}^m G_i \leq \frac{\alpha}{p} m\right) \leq \exp[-300m] \qquad \forall m \geq 1$$

*holds. Moreover, for every $C > 0$ there exists $1 < A_0(C) < \infty$, such that for every $A > A_0$*

$$(58) \qquad P\left(\sum_{i=1}^m G_i > \frac{A}{p} m\right) \leq \exp[-Cm] \qquad \forall m \geq 1.$$

Recall the definition of the event $F$: $\forall v \in V_{1\%}$, $N_v^- \leq 2\alpha_1 N_{1v}$.

Let $u$ be a nonnegative integer. Let us consider an $(-u)$-cell. Recall the random variables $\xi_i$ and $\eta_i$ as in Section 6.1 and recall the random variable $T$ as in (41), which is the smallest time $T$ such that $\xi_i \neq 0$ and $\eta_{u+i} \neq 0$. Let $T_x(j)$ be the index of $j$th $\xi_i$ such that $\xi_i \neq 0$. So

$$T_x(1) = \min\{i \geq 1 : \xi_i \neq 0\}, \qquad \ldots,$$
$$T_x(j+1) = \min\{i > T_x(j) : \xi_i \neq 0\}.$$

Let

$$(59) \qquad \rho^- := \min\{j = 1, 2, \ldots : \eta_{u+T_x(j)} \neq 0\}.$$

Hence $\rho^-$ is the number of $\xi_i$'s in the cell that are not 0 (including the one that is aligned). With this notation,

$$T = T_x(\rho^-).$$

For a $(-u)$-cell, the number of 0-s in $X$ is smaller then the number of 0's in $Y$. Let us estimate (from below) the number of 1's inside the $X$-side, $N_1^-$. This number does not count the aligned one, so the number is clearly at least $\rho^- - 1$, that is, $N_1^- \geq \rho^- - 1$, where the equality holds if and only if

$$\xi_{T_x(j)} = 1, \qquad j = 1, \ldots, \rho^- - 1.$$

The random variable $\rho^-$ has geometric distribution with parameter $\varepsilon$. Indeed, since $X$ and $Y$ are independent, from the right-hand side of (59) follows

$$P(\rho^- = n) = P(\eta_{u+T_x(1)} = 0, \ldots, \eta_{u+T_x(n-1)} = 0, \eta_{u+T_x(n)} \neq 0)$$
$$= (1-\varepsilon)^{n-1}\varepsilon.$$

Let $v = (v_1, \ldots, v_k)$. Recall that $N_v^-$ is the number of ones on the sides with fewer 0's of nonzero cells. At first, we give a lower bound on $N_v^-$.

LEMMA 7.6. *There exists a $\gamma > 0$, not depending on $\varepsilon$ and $\varepsilon_1 > 0$ such that for every $v = (v_1, \ldots, v_k) \in V_{1\%}$ we have*

(60) $$P(F_{1v}^c) \leq \exp[-3k] \qquad \text{where } F_{1v} = \left\{ N_v^- \geq \frac{\gamma}{\varepsilon} k \right\},$$

*provided $\varepsilon < \varepsilon_1$.*

PROOF. Let $v = (v_1, \ldots, v_k) \in V_{1\%}$. Let $I$ be the index set of nonzero cells, $|I| \geq 0.01k$. Let us estimate (below) the number of 1's in the side of fewer 0's:

$$N_v^- = \sum_{i=1}^{|v|} N_v^-(i).$$

For a cell $v_i \neq 0$, we have that $N_v^-(i) \geq \rho_i^- - 1$, where $\rho_i^-$, $i \in I$ are geometrically distributed random variables with parameter $\varepsilon$ as in (59). So,

(61) $$N_v^- \geq \sum_{i \in I} \rho_i^- - |I|.$$

Let $\alpha_o$ be as in Lemma 7.5. It does not depend on $\varepsilon$. Let

$$m := 0.01k, \qquad \gamma := \frac{\alpha_o}{200}, \qquad \varepsilon_1 := \frac{\alpha_o}{2}.$$

Without loss of generality, we may assume that the $v_1, \ldots, v_{|I|}$ are nonzero. Thus, if $\varepsilon \leq \varepsilon_1$, then $100\gamma + \varepsilon \leq \alpha_o$, and from (57) we obtain:

$$P(F_{1v}^c) \leq P\left(\sum_{i \in I}(\rho_i^- - 1) \leq \frac{\gamma}{\varepsilon}k\right)$$

$$\leq P\left(\sum_{i=1}^m (\rho_i^- - 1) \leq \frac{\gamma}{\varepsilon}k\right)$$

$$= P\left(\sum_{i=1}^m \rho_i^- \leq \left(\frac{100\gamma}{\varepsilon} + 1\right)m\right)$$

$$= P\left(\sum_{i=1}^m \rho_i^- \leq \frac{100\gamma + \varepsilon}{\varepsilon}m\right)$$

$$\leq P\left(\sum_{i=1}^m \rho_i^- \leq \frac{\alpha_o}{\varepsilon}m\right)$$

$$\leq \exp[-300(0.01)k]$$

$$= \exp[-3k]. \qquad \square$$

Let
$$F_1(k) := \bigcap_{v \in V_{1\%} \cap V(k)} F_{1v} \quad \text{and} \quad F_1 := \bigcap_{k \geq (0.1\varepsilon^2)n} F_1(k).$$

By (50), (60) and Lemma 7.6, it holds: if $\varepsilon \leq \varepsilon_1$, then

$$P(F_1(k)^c) \leq \sum_{v \in V(k)} P(F_{1v}^c) \leq 16^k \exp[-3k] = \exp[(\ln 16 - 3)k] \leq \exp[-0.2k].$$

Hence

$$P(F_1^c) \leq \sum_{k \geq (0.1\varepsilon^2)n} P(F_1(k)^c) \leq \sum_{k \geq (0.1\varepsilon^2)n} \exp[-0.2k]$$

(62)
$$= C_{1,F} \exp[-0.2(0.1\varepsilon^2)n],$$

where

$$C_{1,F} := (1 - \exp[-0.2])^{-1}.$$

Let $v = (v_1, \ldots, v_k) \in V(k)$ be given. Let $C_v(1), \ldots, C_v(k)$ be the corresponding cells. Let $\rho_j$, $1 \leq j \leq k$ be the number of nonzero $\xi_i$'s in the cell $C_v(j)$. Clearly $\rho_1, \ldots, \rho_k$ are independent. The distribution of $\rho_j$ is geometric with parameter $\varepsilon$, if $v_j \leq 0$. Otherwise, there exists a geometric random variable with parameter $\varepsilon$, say $\rho_j^-$ such that $\rho_j^- \leq \rho_j \leq \rho_j^- + v_j$. Since $v \in V(k)$, $\sum_j |v_j| \leq 2k$. Let us estimate from above the quantity $\rho_v := \sum_{j=1}^k \rho_j$.

LEMMA 7.7. *There exist a constant $B$ not depending on $\varepsilon$ such that for every $v = (v_1, \ldots, v_k) \in V_n$ we have*

$$P(F_{2v}^c) \leq \exp[-(\ln 16 + 1)k], \qquad \text{where } F_{2v} := \left\{\rho_v < \frac{B}{\varepsilon}k\right\}.$$

PROOF. Let $B$ be such that $B - 1 > A_0(\ln 16 + 1)$, let $v = (v_1, \ldots, v_k) \in V(k)$. By (58),

$$P(F_{2v}^c) = P\left(\sum_{j=1}^k \rho_j \geq \frac{B}{\varepsilon}k\right)$$

$$\leq P\left(\sum_{j:v_j \leq 0} \rho_j + \sum_{j:v_j > 0} (\rho_j^- + v_j) \geq \frac{B}{\varepsilon}k\right)$$

$$\leq P\left(\sum_{j:v_j \leq 0} \rho_j + \sum_{j:v_j > 0} \rho_j^- + 2k \geq \frac{B}{\varepsilon}k\right)$$

$$\leq P\left(\sum_{j:v_j \leq 0} \rho_j + \sum_{j:v_j > 0} \rho_j^- \geq \frac{B - 2\varepsilon}{\varepsilon}k\right)$$

$$\leq P\left(\sum_{j:v_j \leq 0} \rho_j + \sum_{j:v_j > 0} \rho_j^- \geq \frac{B - 1}{\varepsilon}k\right)$$

$$\leq \exp[-(\ln 16 + 1)k]. \qquad \square$$

Let

$$F_2(k) := \bigcap_{v \in V(k)} \left\{\rho_v < \frac{B}{\varepsilon}k\right\}, \qquad F_2 := \bigcap_{k \geq (0.1\varepsilon^2)n} F_2(k).$$

Then, similarly to (62),

$$P((F_2(k))^c) \leq \sum_{v \in V(k)} P(F_{2v}^c) \leq \exp[-k((\ln 16 + 1) - \ln 16)] = \exp[-k],$$

$$P(F_2^c) \leq C_{2F} \exp[-0.1\varepsilon^2 n],$$

where

$$C_{2F} := (1 - \exp[-1])^{-1}.$$

Next, using Lemma 7.7, we estimate from above the random number of ones in the $X$-side of the cells $C_v(1), \ldots, C_v(|v|)$.

LEMMA 7.8. *There exists a constant $A < \infty$, independent of $\varepsilon$, such that for every $v = (v_1, \ldots, v_k) \in V_n$, we have*

$$P\left(\sum_{j=1}^{\pi(k)} X_j > \frac{Ak}{\varepsilon(1-\varepsilon)}\right) \leq 2\exp[-(\ln 16 + 1)k].$$

PROOF. Let $v = (v_1, \ldots, v_k) \in V(k)$. Note that
$$P(\xi_i = k | \xi_i \neq 0) = \varepsilon^{k-1}(1-\varepsilon), \qquad k = 1, 2, \ldots.$$
The number of 1's on the $X$-side of the cell $C_v(j)$ is

(63) $$\sum_{i=1}^{\rho(j)-1} G_i + 1,$$

where $G_i$ are geometrically distributed r.v-s with parameter $1 - \varepsilon$ independent of $\rho(j)$. Here $\rho(j) - 1$ is the number of $\xi$'s inside the cell $C_v(j)$ and the additional one is the one that is matched. Hence,

(64) $$\sum_{j=1}^{\pi(k)} X_j = \sum_{i=1}^{\rho_v - k} G_i + k \leq \sum_{i=1}^{\rho_v} G_i.$$

Let $B$ be as in the previous lemma and let $A$ be large enough so that
$$\frac{A}{B} > A_o\left(\frac{(\ln 16 + 1)}{B}\right)$$
and define
$$F_{3v} := \left\{\sum_{i=1}^{B/\varepsilon k} G_i < \frac{A}{\varepsilon(1-\varepsilon)}k\right\}.$$

From Lemma 7.5 with $m = \frac{B}{\varepsilon}k$

$$P(F_{3v}^c) = P\left(\sum_{i=1}^{B/\varepsilon k} G_i \geq \frac{Ak}{\varepsilon(1-\varepsilon)}\right) = P\left(\sum_{i=1}^{B/\varepsilon k} G_i \geq \frac{k}{(1-\varepsilon)}\frac{B}{\varepsilon}\frac{A}{B}\right)$$
$$= P\left(\sum_{i=1}^{m} G_i \geq \frac{mA}{B(1-\varepsilon)}\right) \leq \exp\left[-\frac{(\ln 16 + 1)}{B}m\right]$$
$$< \exp\left[-\frac{(\ln 16 + 1)\varepsilon}{B}m\right] = \exp[-(\ln 16 + 1)k].$$

Due to (64),
$$F_{2,v} \cap F_{3,v} \subset \left\{\sum_{j=1}^{\pi(k)} X_j \leq \frac{Ak}{\varepsilon(1-\varepsilon)}\right\} =: F_{4,v}.$$

Lemma 7.7 finishes the proof. □

Let
$$F_4(k) := \bigcap_{v \in V(k) \cap V_{1\%}} F_{4v}, \qquad F_4 := \bigcap_{k \geq (0.1\varepsilon^2)n} F_4(k).$$

Then, just like in (62),

(65)    $P(F_4^c(k)) \leq 2\exp[-k((\ln 16 + 1) - \ln 16)] = 2\exp[-k],$

(66)    $P(F_4^c) \leq 2C_{2F} \exp[-0.1(\varepsilon^2)n].$

LEMMA 7.9.  *There exists $\alpha_1 > 0$, independent of $\varepsilon$, as well as a constant $C_F < \infty$ such that*
$$P(F^c) \leq C_F \exp[-0.02\varepsilon^2 n],$$
*provided $\varepsilon < \varepsilon_1$, where $\varepsilon_1$ is as in Lemma 7.6.*

PROOF.  Let $\varepsilon < \varepsilon_1$ and $v \in V_{1\%}$. We have
$$F_{1,v} \cap F_{4,v} \subset \left\{ N_v^- \geq \frac{(1-\varepsilon)\gamma}{A} \sum_{j=1}^{\pi(|v|)} X_j \right\}.$$

So,
$$F_1 \cap F_4 = \left( \bigcap_{v \in V_{1\%}} F_{1,v} \right) \cap \left( \bigcap_{v \in V_{1\%}} F_{4,v} \right) = \bigcap_{v \in V_{1\%}} (F_{1,v} \cap F_{4,v})$$
$$\subset \bigcap_{v \in V_{1\%}} \left\{ N_v^- \geq \frac{(1-\varepsilon)\gamma}{A} \sum_{j=1}^{\pi(|v|)} X_j \right\} =: F_x$$

and by (62) and (66)
$$P(F_x^c) \leq P(F_1^c) + P(F_4^c) \leq C_{F1} \exp[-0.02\varepsilon^2 n] + 2C_{F2} \exp[-0.1\varepsilon^2 n].$$

By symmetry, $P(F_y^c) \leq C_{F1} \exp[-0.02\varepsilon^2 n] + 2C_{F2} \exp[-0.1\varepsilon^2 n]$, where
$$F_y := \left\{ N_v^- \geq \frac{(1-\varepsilon)\gamma}{A} \sum_{j=1}^{\nu(|v|)} Y_j \right\}.$$

Thus
$$F_x \cap F_y \subset \left\{ 2N_v^- \geq \frac{(1-\varepsilon)\gamma}{A} \left( \sum_{j=1}^{\pi(|v|)} X_j + \sum_{j=1}^{\nu(|v|)} Y_j \right) \right\} \subset \{N_v^- \geq 2\alpha_1 N_{1v}\} = F,$$

where

(67) $$\alpha_1 := \frac{\gamma}{8A} \leq \frac{(1-\varepsilon)\gamma}{4A},$$

provided $\varepsilon \leq 0.5$ and

$$P(F^c) \leq 2C_{F1}\exp[-0.02\varepsilon^2 n] + 4C_{F2}\exp[-0.1\varepsilon^2 n]$$
$$< (2C_{F1} + 4C_{F2})\exp[-0.02\varepsilon^2 n]. \qquad \square$$

7.4. *The event $G$.* We use the notation introduced in the previous subsection. Let $\alpha_1$ be as in (67). Fix $0 < \alpha_2 < \alpha_1$.

LEMMA 7.10. *There exists a constant $C_G < \infty$ and $\varepsilon_2(\alpha_2) > 0$ such that for every $\varepsilon \leq \varepsilon_2$*

$$P(G^c) \leq C_G \exp[-(300 - \ln 16)(0.1)\varepsilon^2 n].$$

PROOF. Let $v \in V(k)$. From (64)

$$\sum_{j=1}^{\pi(k)} X_j = \sum_{j=1}^{\rho_v - k} G_i + k \geq \rho_v = \sum_{i=1}^{k} \rho_j \geq \sum_{i=1}^{k} \rho_j^-.$$

Let

$$G_v := \left\{ k \leq \alpha_2 \sum_{j=1}^{\pi(k)} X_j \right\}.$$

Then

$$P(G_v^c) \leq P\left( \sum_{i=1}^{k} \rho_j^- < \frac{k}{\alpha_2} \right) = P\left( \sum_{i=1}^{k} \rho_j^- < \frac{\varepsilon}{\alpha_2} \frac{1}{\varepsilon} k \right).$$

Let $\alpha_o$ be as in Lemma 7.5. Let $\varepsilon_2 := \alpha_2 \alpha_o$. Note that $\alpha_2 < 0.5$, so $\varepsilon_2$ is smaller than $\varepsilon_1$ defined in Lemma 7.6. Recall that $\rho_i^-$ are i.i.d. random variables with $G(\varepsilon)$ distribution. Then, by (57), for every $\varepsilon \leq \varepsilon_2$,

$$P(G_v^c) \leq \exp[-300k].$$

Let

$$G(k) := \bigcap_{v \in V(k)} G_v,$$

$$\bigcap_{k \geq 0.1\varepsilon^2} G(k) = \bigcap_{v \in V_n} G_v \subset \bigcap_{v \in V_{1\%}} \left\{ |v| \leq \alpha_2 \sum_{j=1}^{\pi(|v|)} X_j \right\} =: G_x.$$

There exists a constant $0.5C_G$ such that, for $\varepsilon \leq \varepsilon_2$,
$$P(G_v^c(k)) \leq \exp[-(300 - \ln 16)k],$$
$$P(G_x^c) \leq 0.5 C_G \exp[-(300 - \ln 16)(0.1\varepsilon^2)n].$$
Similarly $P(G_y^c) \leq 0.5 C_G \exp[-(300 - \ln 16)(0.1\varepsilon^2)n]$, where
$$G_y := \bigcap_{v \in V_{1\%}} \left\{ |v| \leq \alpha_2 \sum_{j=1}^{\nu(|v|)} Y_j \right\}.$$
Since $G := G_x \cap G_y$, we have that
$$P(F^c) \leq C_G \exp[-(300 - \ln 16)(0.1\varepsilon^2)n],$$
provided $\varepsilon \leq \varepsilon_2$. □

7.5. *The event K.*

LEMMA 7.11. *Assume $\varepsilon < \varepsilon_3$, where $\varepsilon_3$ is as in Lemma 6.2. There exists a constant $C_K$ such that*
$$P(K^c) \leq C_K \exp[-c_K n],$$
*where $c_K > 0$ is a constant, depending on $\varepsilon$.*

PROOF. Let $v$ be an optimal alignment of $X$ and $Y$ and denote by $R_v$ the number of ones after the last cell:
$$R_v := \sum_{i=\pi(|v|)+1}^{n} X_i + \sum_{i=\nu(|v|)+1}^{n} Y_i.$$
Let
$$\beta := (0.1)\varepsilon^2.$$
We shall show that with high probability, there exists an optimal alignment such that all the ones after the last cell are contained in the interval $[n - \beta n + 1, n]$ (without loss of generality assume that $\beta n$ is an integer). In other words, we shall prove that the following event has big probability:
$$(68) \qquad K_1 := \{ \exists v \in V^* : \pi(|v|) \geq n - \beta n, \ \nu(|v|) \geq n - \beta n \},$$
where $V^*$ is the set of optimal alignments.

Suppose $K_1$ holds and let $v \in V^*$ be such that $\pi(|v|) \geq n - \beta n$ and $\nu(|v|) \geq n - \beta n$. Then the number of ones after the last cell of $v$ is clearly at most $2\beta n$, since there are at most $2\beta n$ symbols after the last cell. Thus, $R_v \leq 2\beta n$. Recall that
$$N_{1v} = \sum_{i=1}^{\pi(|v|)} X_i + \sum_{i=1}^{\nu(|v|)} Y_i.$$

Obviously $N_{1v} \geq 2|v|$ and if $E_4$ holds, then every optimal $v$ satisfies $|v| \geq (0.1)\varepsilon^2 n$. Hence, if $K_1 \cap E_4$ holds, then there exists an $v \in V^*$ such that

$$R_v \leq 2\beta n = 2(0.1)\varepsilon^2 \leq 2|v| \leq N_{1v},$$

implying that

$$P(K^c) \leq P(K_1^c) + P(E_4^c).$$

It remains to show that $K_1^c$ has exponentially small probability in $n$. Define

$$K_2 := \left\{ \exists s, t > n - \frac{\beta}{3}n : s + \sum_{i=s+1}^{n} X_i = t + \sum_{i=t+1}^{n} Y_i, \ Y_t = X_s = 1 \right\}.$$

Recall that after the last cell, only the zeros can be aligned. If, in an interval one aligns only the zeros, it can be done in the following manner. Start from the last pair of zeros and align them. Then, disregarding all the ones, take the second last pair of zeros (i.e., the second last zero in $X$ and the second last zero in $Y$) and align them. Then, align the third last pair of zeros and so on. Doing so, the maximum number of zero-pairs (in the given interval) can be obtained. If the event $K_2$ holds, then in the interval $[n - \frac{\beta}{3}n + 1, n]$, the described way of aligning zeros allows to align a pair of ones without disturbing the alignment of zeros. This violates the optimality, hence we immediately have the following implication: if $K_2$ holds, then for any optimal alignment $v$ either $\pi(|v|) \geq n - \frac{\beta}{3}n$ or $\nu(|v|) \geq n - \frac{\beta}{3}n$. Unfortunately, this is not enough, so we define two more events:

$$K_3^x := \left\{ \sum_{s=n-\beta n+1}^{n-2/3\beta n} X_s \geq 1 \right\} \cap \left\{ \frac{2}{3}\beta n - \sum_{s=n-2/3\beta n+1}^{n} X_s \geq \frac{1}{3}\beta n \right\},$$

$$K_3^y := \left\{ \sum_{s=n-\beta n+1}^{n-2/3\beta n} Y_s \geq 1 \right\} \cap \left\{ \frac{2}{3}\beta n - \sum_{s=n-2/3\beta n+1}^{n} Y_s \geq \frac{1}{3}\beta n \right\}.$$

The event $K_3^x$ states that among

$$X_{n-\beta n+1}, \ldots, X_{n-2/3\beta n}$$

there is at least one and, at the same time, among

$$X_{n-2/3\beta n+1}, \ldots, X_n$$

there are at least $\frac{1}{3}\beta n$ zeros. The event $K_3^y$ is symmetric.

Suppose $K_3^x$ holds. Let $v$ be an optimal alignment such that $\nu(|v|) \geq n - \frac{\beta}{3}n$ and $\pi(|v|) < n - \beta n$. Then there exists another optimal alignment $v'$ such that $\nu(|v'|) = \nu(|v|)$ and $\pi(|v'|) \geq n - \beta n$. Indeed, since $\nu(|v|) \geq n - \frac{\beta}{3}n$, the number of aligned 0's after last cell is at most $\frac{\beta}{3}n$ [the maximal number

of 0's on $Y$-side after $\nu(|v|)$]. By $K_3^x$, we can align all those 0's from the $Y$-side with the zeros on $X$-side that lie on $X_{n-2/3\beta n+1}, \ldots, X_n$. After such a realignment, the situation is the following: as previously, $Y_{\nu(|v|)}$ is aligned with $X_{\pi(|v|)}$. However, all the zeros after $Y_{\nu(|v|)}$ (on $Y$-side) are aligned with the zeros after $X_{n-2/3\beta n+1}$ (on $X$-side). The score remains optimal. By $K_3^x$, again, among $X_{n-\beta n+1}, \ldots, X_{n-2/3\beta n}$, there is at least one 1. Thus, without changing the score, we can align $Y_{\nu(|v|)}$ with this 1. Since the location of this 1 is at least $X_{n-\beta n+1}$, we now have a new optimal alignment $v'$ such that $\nu(|v'|) = \nu(|v|)$ and $\pi(|v'|) \geq n - \beta n$.

We have proven that
$$K_2 \cap K_3^x \cap K_3^y \subset K_1.$$

It remains to prove that $K_2$, $K_3^x$ and $K_3^y$ hold with big probability.

Clearly,
$$P(K_3^{xc}) \leq P\left(\sum_{s=n-\beta n+1}^{n-2/3\beta n} X_s = 0\right) + P\left(\sum_{s=n-2/3\beta n+1}^{n} X_s > \frac{1}{3}\beta n\right).$$

The first probability of the right-hand side equals to $\exp[\ln(1-\varepsilon)\frac{1}{3}\beta n]$; by Höffding's inequality, the second probability is bounded by $\exp[-\frac{4}{3}(\frac{1}{2} - \varepsilon)^2 \beta n]$. Thus,
$$P((K_3^x \cap K_3^y)^c) \leq 2\exp[\ln(1-\varepsilon)\frac{1}{3}\beta n] + 2\exp[-\frac{4}{3}(\frac{1}{2} - \varepsilon)^2 \beta n].$$

Finally, let us bound $P(K_2^c)$. The event $K_2$ essentially states that among i.i.d. Bernoulli $B(1,\varepsilon)$ random variables
$$X_1, \ldots, X_{\beta/3n}, Y_1, \ldots, Y_{\beta/3n}$$
after aligning all 0's, one can align an additional pair of ones. In terms of $\xi_i$'s and $\eta$'s, $K_2 = K_2^x \cap K_2^y$, where
$$K_2^x := \left\{\sum_{j=1}^{T-1} \xi_j + T \leq \frac{\beta}{3}n\right\}, \qquad K_2^y := \left\{\sum_{j=1}^{T-1} \eta_j + T \leq \frac{\beta}{3}n\right\}.$$

Recall that $T$ is the stopping time that shows the first time a pair of ones between the 0's occurs. Then $T-1$ is the number of 0's before the first alignment of ones and $\sum_{j=1}^{T-1} \xi_j$ is the number of $X_i$'s before the first alignment of ones. If their sum is smaller than $\frac{\beta}{3}n$, then the $X$-part of the aligned pair occurs before $\frac{\beta}{3}n$. The event $K_2^y$ is analogous.

To bound the events $K_2^x$ and $K_2^y$, we use Lemma 7.5. Let $A_0(1)$ be as in Lemma 7.5 and define
$$\delta := \frac{\beta}{3}\frac{(1-\varepsilon)}{A_0(1)} < \frac{\beta}{3}.$$

Clearly

$$P(K_2^{xc}) \le P(T > \delta n) + P\left(\sum_{i=1}^{\delta n} \xi_i + \delta n > \frac{\beta}{3}n\right).$$

Since $T$ is a geometric random variable with parameter $\varepsilon^2$, $P(T > \delta n) \le \exp[\ln(1-\varepsilon^2)\delta n]$. Since $G_i := \xi_i + 1$ are geometric random variables with parameter $(1-\varepsilon)$, by (58), we have

$$P\left(\sum_{i=1}^{\delta n}\xi_i + \delta n > \frac{\beta}{3}n\right) = P\left(\sum_{i=1}^{\delta n} G_i > \frac{\beta}{3}n\right) = P\left(\sum_{i=1}^{\delta n} G_i > \frac{A_0(1)}{1-\varepsilon}\delta n\right)$$
$$\le \exp[-\delta n].$$

To sum up:

$$P(K^c) \le P(K_2^{xc}) + P(K_2^{yc}) + P(K_3^{xc}) + P(K_3^{yc}) + P(E_4^c)$$
$$\le 2\exp[\ln(1-\varepsilon)\tfrac{1}{3}\beta n] + 2\exp[-\tfrac{4}{3}(\tfrac{1}{2}-\varepsilon)^2\beta n]$$
$$+ 2\exp[\ln(1-\varepsilon^2)\delta n] + 2\exp[-\delta n]$$
$$+ 5\exp[-an] + 4\exp[-2(0.8\varepsilon)^2\varepsilon n],$$

since $\varepsilon < \varepsilon_3$, so $P(E_4^c) \le 5\exp[-an] + 4\exp[-2(0.8\varepsilon)^2\varepsilon n]$ by Lemma 6.2. □

## APPENDIX

PROPOSITION A.1. *Let $G_1, \ldots, G_m$ be i.i.d. geometrically distributed random variables with parameter $p$. Then for every $A > 1$ and $\alpha < 1$, there exists $C(A) := A - 1 - \log A$ and $C(\alpha) := \alpha - 1 - \log \alpha$ such that such that*

$$(69) \qquad P\left(\sum_{i=1}^m G_i > \frac{A}{p}m\right) \le \exp[-C(A)m],$$

$$(70) \qquad P\left(\sum_{i=1}^m G_i \le \frac{\alpha}{p}m\right) \le \exp[-C(\alpha)m].$$

PROOF. Let us recall (55). Let $A > 1$, $n = \frac{A}{p}m$ and $a = \frac{p}{A} < p$. From (55), we get

$$P\left(\sum_{i=1}^m G_i > \frac{A}{p}m\right) = P\left(\sum_{j=1}^{A/pm} Y_j < m\right) = P\left(\sum_{j=1}^n Y_j < an\right)$$
$$\le \left(\frac{p}{a}\right)^{an} \exp[(a-p)n] = A^m \exp[(1-A)m]$$
$$= \exp[(\ln A - (1-A))m],$$

where $Y_i$ are i.i.d. Bernoulli random variables with parameter $p$. This finishes the proof of (69).

If $p > \alpha$, then (70) trivially holds. If $p = \alpha$, then the probability in (70) equals $p^m = \exp[(\ln \alpha)m] = \exp[-\ln \frac{1}{\alpha} m]$. Hence, we consider only the case $p < \alpha < 1$. From (55), it easily follows: let $X_i \sim B(1, p)$. Then, with $1 > a \geq p$,

$$(71) \quad P\bigg(\sum_{i=1}^{n} X_i \geq na\bigg) \leq \exp\bigg[\bigg(a \ln\bigg(\frac{p}{a}\bigg) + (a - p)\bigg)n\bigg].$$

We have that

$$(72) \quad P\bigg(\sum_{i=1}^{m} G_i \leq \frac{\alpha}{p} m\bigg) = P\bigg(\sum_{j=1}^{\alpha/pm} Y_j \geq m\bigg) = P\bigg(\sum_{j=1}^{\alpha/pm} Y_j \geq \frac{\alpha}{p} m \frac{p}{\alpha}\bigg).$$

With $n := \frac{\alpha}{p} m$ and $a := \frac{p}{\alpha}$, the inequality (71) states

$$P\bigg(\sum_{i=1}^{m} G_i \leq \frac{\alpha}{p} m\bigg) \leq \exp\bigg[\bigg(\frac{p}{\alpha} \ln \alpha + \bigg(\frac{p(1 - \alpha)}{\alpha}\bigg)\bigg)n\bigg]$$
$$= \exp[(\ln \alpha + 1 - \alpha)m]. \quad \square$$

PROOF OF LEMMA 7.5. The right-hand side of (70) is smaller than $\exp[-300]$, provided

$$\alpha \leq \exp[-301] =: \alpha_o.$$

That proves (57). To get (58), note that for every $C > 0$, it is possible to choose $A$ so big that $\ln A - (1 - A) < -C$. $\square$

**Acknowledgment.** The authors would like to thank the anonymous referee for valuable suggestions and remarks.

INSTITUTE OF MATHEMATICAL STATISTICS
UNIVERSITY OF TARTU
J. LIIVI 2-513
50409, TARTU
ESTONIA
E-MAIL: jyril@ut.ee

SCHOOL OF MATHEMATICS
GEORGIA TECH
ATLANTA, GEORGIA 30332-0160
USA
E-MAIL: matzing@math.gatech.edu